\documentclass{article}
\usepackage{amsfonts}
\usepackage{amssymb,amsmath,multirow}
\usepackage[dvips]{graphicx}
\usepackage[dvips]{color}
\usepackage{enumitem}
\usepackage{lscape,makecell,rotating,multirow,threeparttable}

\newtheorem{thm}{Theorem}
\newtheorem{lem}{Lemma}
\newtheorem{cor}{Corollary}
\newtheorem{rem}{Remark}
\newtheorem{algo}{Algorithm}
\def\pn{\par\smallskip\noindent}
\def\proof{\pn {Proof.} }
\def\endproof{\hfill \quad{$\Box$}\smallskip}

\title{Optimality Conditions for  Cardinality-Constrained  Programs  and a SCA Method \footnotemark[2]}
\author{Zhongyi Jiang$^{*b,a}$\ \ \  Baiyi Wu$^c$ \ \ \ Qiying Hu$^{a}$\\
\small( a. School of Management, Fudan University, Shanghai, P.R.China.)\\
\small( b. School of Information Science, Changzhou University,
        Changzhou, P.R.China)\\
\small(c. School of Finance, Guangdong University of Foreign Studies, Guangzhou, P.R.China.)
}

\begin{document}
\normalsize
\date{}
\maketitle
\renewcommand{\thefootnote}{\fnsymbol{footnote}}
\footnotetext[2]{This research is supported in part by  National Natural Science Foundation of China under grant 11371103.}
\footnotetext[1]{Corresponding author.}
\footnotetext[0]{E-mail addresses:
jiangzy\symbol{64}fudan.edu.cn (Z.Y. Jiang), baiyiwu\symbol{64}outlook.com (B.Y. Wu), qyhu\symbol{64}fudan.edu.cn (Q.Y. Hu) }

\vskip 0.3cm \noindent {\bf Abstract:}
We study  a cardinality-constrained optimization problem with nonnegative variables in this paper. This problem is often encountered in practice.
 Firstly we study some  properties on the optimal solutions of this  optimization problem under some conditions. An equivalent reformulation of the problem under consideration is proposed. Based on the reformulation, we present a successive convex approximation
method for the cardinality constrained optimization problem.
   We prove that the method
converges to a KKT point of the reformulation problem. Under some conditions, the KKT points of the reformulation problem are local optimizers
of the original problem. Our  numerical results on a limited diversified mean-variance
portfolio selection problem demonstrate some promising results.





\section{Introduction}
 The recent literature has witnessed an increasing attention to optimization problems with a sparsity constraint (where the cardinality of the  decision vector is bounded from above), due to their wide spectra of applications in, for example, portfolio selection \cite{bienstock96,chang00,gaoli2013,shaw08,woodside2011,zheng2014},
subset selection in multivariate regression \cite{arthanari93,miller02}, signal processing
and compressed sensing \cite{bruckstein2009}. Given $x$ $\in$ $\Re^n$, the $\ell_0$-(quasi) norm $\|x\|_0$ denotes the number of nonzero components of $x$ and the sparsity constraint $\|x\|_0\le K$ with $1\le K<n$
   is also called the {\em cardinality} constraint. A related, albeit different, optimization problem is to find the sparsest solutions of linear systems $Ax=b$, i.e.  minimizing $\|x\|_0$ subject to $Ax=b$ (see, e.g., \cite{jokar2008} and the references therein). Because of the non-tractability of the so-called zero norm $\|x\|_0$, researchers, see for example \cite{CRT06a, CRT06b,ZhaoLi2012}, have proposed to use $\ell_1$ norm to develop good approximate algorithms.

We study in this paper the following cardinality-constrained optimization problem with nonnegative variables:
   \begin{align*}
     {\rm(P)}~~~~
     \min~~
     &f(x)\\
     {\rm s.t.}~~~
     &g(x)\le 0,\\
     &\|x\|_0\le K,\\
     &x\geq 0,
    \end{align*}
   where $f:\Re^n\to \Re$ is a differentiable convex function, and
   $g(x)=(g_1(x),\ldots,g_m(x))^T$
   with all $g_i:\Re^n\to \Re$ being  differentiable convex functions.
    By introducing a 0-1 variable $y_i$ for each $x_i$,
   the cardinality constraint can be represented by equivalent mixed-integer constraints
   $\sum_{i=1}^n y_i\le K$ and $l_iy_i\le x_i\le u_iy_i$, $y_i\in\{0,1\}$, $i=1,\ldots,n$,
   where $l_i$ and $u_i$ are lower and upper bounds of $x_i$, respectively.
   Therefore, problem (P) can be reformulated as a mixed-integer convex program.

An important subclass of problem (P) is the following cardinality-constrained quadratic program:
  \begin{align*}
     {\rm(QP)}~~~~
     \min~~
     &x^TQx+c^Tx\\
     {\rm s.t.~~~}
     &Ax\le b,\\
     &\|x\|_0\le K,
   \end{align*}
where $Q$ is an $n\times n$ symmetric positive semidefinite matrix.
The cardinality-constrained portfolio selection model is essentially a special case of
(QP) where the cardinality constraint confines the total number of different assets
in the optimal portfolio. Problem ${\rm(QP)}$ can also be reformulated as a mixed-integer convex program.   Using different
relaxations and bounding techniques, various branch-and-bound methods have been proposed for solving
the mixed-integer quadratic program reformation of (QP)
 (see, e.g., \cite{bertsimas09,bienstock96,bonami09,li06,shaw08,vielma08,gaoli2013}).
A mixed-integer quadratically constrained quadratic program reformulation
is derived in \cite{cui2012} for a class of cardinality-constrained portfolio selection
problems where the asset returns are driven by factor models.

Index tracking in passive portfolio management is one of the important applications of the quadratic model (QP), in which
a small set of assets is selected to track the performance of the market benchmark index. Different  approximation methods for solving the index tracking problem can be found in
\cite{jansen2002,coleman2006,zheng2012}.

Problem (P) has been proved to be  NP-hard in \cite{bienstock96}.  Notice  that testing
the feasibility of (QP) is already NP-complete when $X=\{x\mid Ax\leq b\}$ and $A$ has three rows \cite{bienstock96}.
 So it is a very tough job to find a feasible solution of the problem, let alone  a local optimal solution.

A penalty form of problem (P), which is often called
{\em regularization}
formulation of (P), is  used to  find sparse solutions  by
attaching the $\ell_0$ function in the objective function as follows:
\begin{align*}
{\rm(P_\mu)}~~~~
\min~~
&f(x)+\mu\|x\|_0\\
{\rm s.t.}~~~
&g(x)\le 0,
\end{align*}
where $\mu>0$ is a regularization parameter.
It is  shown in \cite{jansen2002} that for each $1\le K<n$, the
optimal solution of (P) can be generated via solving ${\rm(P_\mu)}$ for some $\mu>0$
under some mild conditions.
We can solve ${\rm(P_\mu)}$ (approximately) with different values of $\mu$
and eventually find an optimal solution of (P) which satisfies $\|x\|_0\leq K$.
The $l_0$ norm $\|x\|_0$ can be approximated by the $\ell_p$ norm,
   $\|x\|^p_p$, with $p\in(0,1)$, which leads to an $\ell_2$-$\ell_p$ minimization problem when the objective function is quadratic.
   The lower bound theory of nonzero entries of $\ell_2$-$\ell_p$
   minimization is discussed in \cite{chen2011,chen2010a,chen2010b}.
   An interior-point potential reduction algorithm
  is proposed   in \cite{ge2011} to search for a local solution of $\ell_2$-$\ell_p$ minimization.
A DC (difference of convex functions) approximation method is proposed in  \cite{zheng2012} to find approximation solutions of
${\rm(P_\mu)}$.

Several different optimality conditions and  algorithms for nonlinear optimization problems with sparsity constraints are presented in \cite{Beck2013,Lu2013}.
Recently,  some different complementarity-type reformulations and regularization methods for the optimization problems with sparsity constraints are presented in \cite{Feng2013, Burdakov2014}. A DC (difference of
convex) approaches for linear programs with complementarity constraints was presented in \cite{Francisco2018}. The DC method can also be used to solve sparsity constraints optimization problems by transforming the sparsity constraints into complementarity constraints.

In this paper, we propose a  reformulation for problem (P) and then develop  a successive convex approximation method by  successively
linearizing the  term in the reformulation. Under some constraint qualifications, we prove necessary optimal conditions of our reformulation and establish the convergence of
the sequence of approximate solutions to a KKT point of
problem (P), which is also a local optimal solution of problem (P). We test our method on the limited diversified mean-variance
portfolio selection problem. Our preliminary computational results do demonstrate some promising properties of our proposed solution scheme.

Our main contributions include the following two aspects: Firstly, we derive conditions under which the cardinality constraint is binding at the optimal point; Secondly we give a transformation of the cardinality constraint. Based on this transformation, we give a successive convex approximation method to the cardinality constrained optimization problems. Compared with the regularized method, our method can get a local optimal solution which exactly satisfies the cardinality constraint under certain conditions.

The rest of the paper is organized as follows. In Section 2,  we derive
some properties for the optimal solution of the cardinality constrained optimization problem (P). In Section 3, we propose an equivalent reformulation of problem (P) and investigate both the global and local optimality conditions.  In Section 4, we  develop a successive convex approximation method based on our reformulation and establish then some convergence results of the method.
After reporting computational results in Section 5 for the  limited diversified mean-variance
portfolio selection problem,
we give some concluding remarks in Section 6.

\section{Properties of the optimal solution }
In this section, we derive some properties of the optimal solutions to problem (P). We will prove these properties without the  nonnegative constraints. It is obvious that these properties are still true with the nonnegative constraints.  Firstly, we would like to derive conditions under which the cardinality constraint is binding at the optimal point. Let ($x,y$) represent point  ($x^T,y^T)^T$ for the convenience in the whole paper.
\begin{thm}\label{th1}
Suppose that $f: \Re^n\to \Re$  and $g_i: \Re^n\to \Re$,  $i=1, \ldots, m$, are all convex functions, and that the inner of the  set $X=\{x\mid g(x)\le 0, \|x\|_0\leq K\}$ is nonempty, i.e. there exists $\bar{x}$ with  $\|\bar{x}\|_0\leq K$
such that $g(\bar{x})>0$.  Let $x^*$ be a local optimal solution of problem $\min\{f(x)\mid g(x)\leq 0,\|x\|_0\leq K\}$.
If $\|{x}'\|_0>K$ for all ${x}'\in \arg \min\{f(x)\mid g(x)\leq 0,\}$,  then there must be $\|{x}^*\|_0=K$.
\end{thm}
\proof
Let $I({x}^*)\subseteq \{1,\ldots,n\}$ be an index set for $i$ satisfying  ${x}_i^*=0$. Then
$x^*$ must be the global optimal solution of the convex optimization problem $\min\{f(x)\mid g(x)\leq 0, x_i=0~ {\rm for}~ i\in I(x^*)\}$, i.e. every local optimal solution of problem $\min\{f(x)\mid g(x)\leq 0,\|x\|_0\leq K\}$ corresponds to a convex optimization problem. We will prove the property by  contradiction through a global optimal solution of $\min\{f(x):g(x)\leq 0\}$.

Suppose $\min \|{x}'\|_0=L$ for all ${x}'\in \arg \min\{f(x):g(x)\leq 0\}$. Then $L>K$ according to our
  assumption. Let $\hat{x}^*\in \arg \min\{f(x)\mid g(x)\leq 0\}$,  $\|\hat{x}^*\|_0=L$,  and $I(\hat{x}^*)\subseteq \{1,\ldots,n\}$ be an index set such that $\hat{x}_i^*=0$ for all $i\in {I(\hat{x}^*)}$. Without loss of generality,
   we  assume that indices 1 and 2 are not in the set $I(\hat{x}^*)$. It is obvious that $\hat{x}^*$ is also the global optimal
   solution of the following problem.
  \begin{align*}
       {\rm(P_0)}~~~
       \min~~
       &f(x)\\
       {\rm s.t.}~~~
       &g(x)\leq0,\\
       &x_i=0,~~i\in I(\hat{x}^*),
  \end{align*}
 Now, we  prove that if $x^*$ is a local optimal solution of the following problem ${\rm (P_1)}$, then $\|x^*\|_0=L-1$.
 \begin{align*}
      {\rm(P_1)}~~~~
      \min~~
      &f(x)\\
      {\rm s.t.}~~~
      &g(x)\leq0,\\
      &x_i=0,~~i\in I(\hat{x}^*),\\
      &\|x\|_0\le L-1.
 \end{align*}
  This will prove the
 property. By contradiction,
 assume  $\|x^*\|_0<L-1$.  Without loss of generality, let $x^*=(0,0,x^*_3,\ldots,x^*_n)$. Under this assumption, we
 will prove that  $x^*$ is also the optimal solution of problem $ {\rm(P_0)}$, which contradicts
 $\|\hat{x}^*\|_0\geq L$. The proof consists of two steps.

 Firstly, we prove  that $x^*$ is also a global optimal solution of the following two problems ${\rm (P_2)}$ and ${\rm (P_3)}$:
 \begin{align*}
      {\rm(P_2)}~~~~
      \min~~&f(x)\\
      {\rm s.t.}~~~
      &g(x)\leq0,\\
      &x_i=0,~~i\in I(\hat{x}^*),\\
      &x_1=0;
 \end{align*}
 \begin{align*}
      {\rm(P_3)}~~~~
      \min~~
      &f(x)\\
      {\rm s.t.}~~~
      &g(x)\leq0,\\
      &x_i=0,~~i\in I(\hat{x}^*),\\
      &x_2=0.
 \end{align*}
By contradiction, suppose $x^*$ is not a global optimal solution of problem ${\rm(P_2)}$, and the global optimal solution of problem ${\rm(P_2)}$ is $\tilde{x}^*$. Without loss of
 generality, let  $\tilde{x}^*=(0,\tilde{x}^*_2,\tilde{x}^*_3,\tilde{x}_4^*\ldots,\tilde{x}^*_n)$. There must be $f(\tilde{x}^*)<f(x^*)$.
 Let $x'=\lambda x^* +(1-\lambda)\tilde{x}^*$, ($0\leq\lambda\leq1$). Then $x'$ is a feasible solution of problem ${\rm(P_1)}$, $ \|x'\|_0\leq L-1$,
 and
 \begin{eqnarray}
 f((\lambda x^* +(1-\lambda)\tilde{x}^*))\leq \lambda f(x^*)+(1-\lambda)f(\tilde{x}^*)<f(x^*).\label{prop1}
  \end{eqnarray}
   Let $\lambda$ get close to 1. Then (\ref{prop1}) is contrary with the fact that $x^*$ is a local optimal solution of problem ${\rm(P_1)}$.
  So $x^*$ is the global optimal solution of problem ${\rm(P_2)}$. By the same way we can  prove that $x^*$ is also the global optimal solution of problem ${\rm(P_3)}$.

 Secondly, we  prove that  ${x}^*$ is a KKT point of   problem ${\rm(P_0)}$, and $x^*$ is also the optimal solution of problem $ {\rm(P_0)}$, which contradicts
 $\|\hat{x}^*\|_0\geq L$.  In fact, the Lagrangian dual  of problem (${\rm P_2}$) is
 \begin{align*}
      {\rm(P_4)}~~~~
      \max~~
      &\theta (u)\\
      {\rm s.t.}~~~
      &u\geq 0,
 \end{align*}
 where $\theta (u)=\inf \{f(x)+\sum_{i=1}^mu_ig_i(x):x_1=0,~x_i=0,{\rm ~for~} i\in I(\hat{x}^*)\}$ with $u_i$ being the $i$th Lagrangian multiplier, $i$ = 1, $\ldots$, $m$. There exist $\bar{x}$ such that $\bar{x}_1=0,~\bar{x}_i=0,{\rm ~for~} i\in I(\hat{x}^*)\}$, and $g(\bar{x})>0$ because the inner of set $X=\{x\mid g(x)\le 0,\|x\|_0\leq k\}$ is nonempty. Then based on Proposition 5.3.1 in \cite{Bertsekas09} (Convex programming duality), we have
 \begin{align*}
    \inf \{f(x):g(x)\leq 0,~x_i=0 {\rm ~for~} i\in I(\hat{x}^*),~x_1=0\}=\sup\{\theta (u): u\geq 0\}.
 \end{align*}
 It is obvious that
 \begin{eqnarray*}
    &&\inf \{f(x):g(x)\leq 0,~x_i=0 {\rm ~for~} i\in I(\hat{x}^*),~ x_2=0\}\\
   &=& \inf \{f(x):g(x)\leq 0,~x_i=0 {\rm ~for~} i\in I(\hat{x}^*),~ x_1=0\}\\
   &=& \sup\{\theta (u): u\geq 0\}.
 \end{eqnarray*}
 Furthermore, the supremum above can be achieved, say at $\bar u$, for the infimum above is finite. Then $f(x^*)=\theta (\bar u)$ and $\bar u^T g(x^*)=0$. Moreover, problem ${\rm (P_2)}$ is equivalent to the following problem
 \begin{align*}
      {\rm(P_5)}~~~~
      \min~~
      &f(x)+\bar u^T g(x)\\
      {\rm s.t.}~~~
      &x_1=0,\\
      &x_i=0,~~i\in I(\hat{x}^*).
 \end{align*}
 Now we  prove that $x^*=(0,0,x^*_3,\ldots,x^*_n)$ is also the optimal solution of the following problem:
 \begin{align*}
      {\rm(P_6)}~~~~
      \min~~
      &f(x)+\bar u^T g(x)\\
      {\rm s.t.}~~~
      &x_2=0,\\
      &x_i=0,~~i\in I(\hat{x}^*).
 \end{align*}
Note that $x^*=(0,0,x^*_3,\ldots,x^*_n)$ is  feasible to problem ${\rm(P_6)}$. By contradiction, if there exists
  $\tilde{x}\in \arg \min\{f(x)+\bar u^T g(x):x_2=0,~x_i=0 {\rm ~for~} i\in I(\hat{x}^*)\}$ such that
 \begin{eqnarray*}
 f(\tilde{x})+\bar u^T g(\tilde{x}) < f(x^*)+\bar u^T g(x^*),
 \end{eqnarray*}
 then there must be
 \begin{eqnarray*}
 f(\tilde{x})+\bar u^T g(\tilde{x}) < f(x^*)+\bar u^T g(x^*)=\theta (\bar u),
 \end{eqnarray*}
 which is contrary to
  \begin{eqnarray*}
 \theta (\bar u)&=&\inf \{f(x)+\sum_{i=1}^m \bar u_ig_i(x):x_1=0,~x_i=0 {\rm ~for~} i\in I(\hat{x}^*)\}\\
 &=&\inf \{f(x)+\sum_{i=1}^m \bar u_ig_i(x):x_2=0,~x_i=0 {\rm ~for~} i\in I(\hat{x}^*)\}.
 \end{eqnarray*}
 Thus,
 $x^*=(0,0,x^*_3,\ldots,x^*_n)$ is also an optimal solution of  problem ${\rm(P_6)}$, and  so problem ${\rm(P_5)}$ is equivalent to problem ${\rm (P_6)}$.\\
 \hspace*{6mm}Note that both  problems ${\rm(P_5)}$ and ${\rm(P_6)}$ are  convex. As $x^*$ is the optimal solution of both ${\rm(P_5)}$ and ${\rm(P_6)}$, $x^*$ is their KKT point as well. The KKT condition of  problem ${\rm(P_5)}$ is
 \begin{eqnarray*}
 \nabla f({x}^*)+\sum_{i=1}^m\bar u_i \nabla g_i({x}^*)+\lambda e_1+\sum_{i\in I(\hat{x}^*)}\lambda'_i e_i=0,
 \end{eqnarray*}
 which can be expressed more specifically as follows,
 \begin{eqnarray}
 &&\frac{\partial f({x}^*)}{\partial x_1} +\sum_{i=1}^m\bar u_i \frac{\partial g_i({x}^*)}{\partial x_1}+\lambda=0,\label{th2}\\
 &&\frac{\partial f({x}^*)}{\partial x_2} +\sum_{i=1}^m\bar u_i \frac{\partial g_i({x}^*)}{\partial x_2} =0,\label{th3}\\
 &&\frac{\partial f({x}^*)}{\partial x_j} +\sum_{i=1}^m\bar u_i \frac{\partial g_i({x}^*)}{\partial x_j}+\lambda'_j  =0, {~\rm for~} j\in I(\hat{x}^*),\label{th4}\\
 &&\frac{\partial f({x}^*)}{\partial x_j} +\sum_{i=1}^m\bar u_i \frac{\partial g_i({x}^*)}{\partial x_j}  =0,
 {~\rm for~} j\not\in I(\hat{x}^*),~ j\neq 1,{\rm ~and~} j\neq 2.  \label{th5}                                                                                          \end{eqnarray}
 On the other hand, the KKT condition of problem ${\rm(P_6)}$ is
 \begin{eqnarray*}
 \nabla f({x}^*)+\sum_{i=1}^m\bar u_i\nabla g_i({x}^*)+\mu e_2+\sum_{i\in I(\hat{x}^*)}\mu'_i e_i=0, \end{eqnarray*}
 which can be expressed more specifically as follows:
 \begin{eqnarray}
 &&\frac{\partial f({x}^*)}{\partial x_1} +\sum_{i=1}^m\bar u_i \frac{\partial g_i({x}^*)}{\partial x_1}=0,\label{th6}\\
 &&\frac{\partial f({x}^*)}{\partial x_2} +\sum_{i=1}^m\bar u_i \frac{\partial g_i({x}^*)}{\partial x_2}+\mu=0,\label{th7}\\
 &&\frac{\partial f({x}^*)}{\partial x_j} +\sum_{i=1}^m\bar u_i \frac{\partial g_i({x}^*)}{\partial x_j}+\mu'_j =0,{~\rm for~} j\in I(\hat{x}^*),\label{th8}\\
 &&\frac{\partial f({x}^*)}{\partial x_j} +\sum_{i=1}^m\bar u_i \frac{\partial g_i({x}^*)}{\partial x_j} =0  {~\rm for~} j\not\in I(\hat{x}^*),~ j\neq 1,{\rm ~and~} j\neq 2.
 \end{eqnarray}
 Comparing (\ref{th2}), (\ref{th3}) and (\ref{th6}), (\ref{th7}) yields $\lambda =0$ and $\mu=0$. Then the KKT condition becomes
  \begin{eqnarray}
  &&\frac{\partial f({x}^*)}{\partial x_1} +\sum_{i=1}^m\bar u_i \frac{\partial g_i({x}^*)}{\partial x_1}=0,\\
  &&\frac{\partial f({x}^*)}{\partial x_2} +\sum_{i=1}^m\bar u_i \frac{\partial g_i({x}^*)}{\partial x_2}=0,\\
  &&\frac{\partial f({x}^*)}{\partial x_j} +\sum_{i=1}^m\bar u_i \frac{\partial g_i({x}^*)}{\partial x_j}+\mu'_j =0,{~\rm for~} j\in I(\hat{x}^*),\\
   &&\frac{\partial f({x}^*)}{\partial x_j} +\sum_{i=1}^m\bar u_i \frac{\partial g_i({x}^*)}{\partial x_j} =0  {~\rm for~} j\not\in I(\hat{x}^*),~ j\neq 1,{\rm ~and~} j\neq 2.
  \end{eqnarray}
Note that
 \begin{eqnarray}
  &&\bar u^T  g({x}^*)=0,\\
  && \mu'_j{x_j}^*=0 {~\rm for~} j\in I(\hat{x}^*),\\
  &&u_i\geq 0, {~\rm for~}i=1,\ldots,n.
  \end{eqnarray}
 We can now conclude  that ${x}^*$ is a KKT point of  problem ${\rm(P_0)}$.
 Note that the equality constraints of problem $ {\rm(P_0)}$ are all linear functions. Then, according to  the KKT necessary conditions \cite{bertsekas2003},
 $x^*$ is also the optimal solution of problem $ {\rm(P_0)}$, which contradicts
 $\|\hat{x}^*\|_0\geq L$.
\endproof \\

Based on Theorem \ref{th1} we can get the following Corollary.
\begin{cor}\label{cor_th1}
Suppose that $f: \Re^n\to \Re$  and $g_i: \Re^n\to \Re$,  $i=1, \ldots, m$, are all convex functions.
  Let $x^L$  be a local optimal solution of problem  $\min\{f(x)\mid g(x)\leq 0,\|x\|_0\leq L\}$,
$x^S$  be a local optimal solution of problem  $\min\{f(x)\mid g(x)\leq 0,\|x\|_0\leq S\}$. If   the inner of the  set $X=\{x\mid g(x)\le 0, \|x\|_0\leq k\}$ is nonempty, $L>S$, and $\|{x}'\|_0>S$ for all ${x}'\in
\arg \min\{f(x):g(x)\leq 0\}$ ,  then  $f(x^L)<f(x^S)$.
\end{cor}
\proof
By contradiction, suppose $f(x^L)\geq f(x^S)$. It is obvious that $x^S$ is  feasible to the following  problem
\begin{align*}
     ~~~~\min~~
     &f(x)\\
     {\rm s.t.}~~~
     &g(x)\leq0,\\
     &\|x\|_0\leq L.
\end{align*}
Then $f(x^L)\geq f(x^S)$, which is contrary to Theorem \ref{th1} because $\|x^S\|_0\leq S<L$.
\endproof \\

Next, we  discuss the properties of local optimal solutions of problem ${\rm (P)}$.
Define $I(x)=\{i\in\{1,\ldots,n\}\mid x_i=0\}$.
For any $z\in\Re^n$, define the following problem parameterized by $z$:
\begin{align*}
{\rm(P}_{I(z)})~~~~
\min~~
&f(x)\\
{\rm s.t.}~~~
&g(x)\leq0,\\
&x_i=0~~{\rm for}~~i\in I(z).
\end{align*}
Suppose that $x^*$ is a local optimal solution of problem ${\rm (P)}$.
Then it is obvious that $x^*$ is also a local optimal solution
of problem ${\rm(P}_{I(x^*)})$.
 If $f(x)$ and $g_i(x)$ ($i=1,\ldots,m$) are all  convex functions, then
${\rm(P}_{I(x^*)})$ is  a convex problem,
 which means that the local optimal solution $x^*$ of  ${\rm(P}_{I(x^*)})$ is also globally optimal.

We show in the following theorem that the converse is also true under some condition.
\begin{thm}\label{local_thm2}
If $x^*$ is a global optimal solution
 of problem ${\rm(P}_{I(x^*)})$ with $\|x^*\|_0= K$, then $x^*$ is also a local optimal solution of problem $\min\{f(x)\mid g(x)\leq 0,\|x\|_0\leq K\}$.
 \end{thm}
 \proof
 It is obvious that $x^*$ is a feasible solution of problem $\min\{f(x)\mid g(x)\leq 0,\|x\|_0\leq S\}$.
Define   $\varepsilon_0=\min\{|x_i^*|\mid x_i^*\neq 0 \}$.
Let  $x'$ be a  feasible solution of problem  ${\rm (P)}$ such that $I(x')\neq I(x^*)$.
 It is easy to see that
    \begin{eqnarray}\label{local_eps}
 \|x'-x^*\|>\varepsilon_0.
   \end{eqnarray}
  Let
  \begin{eqnarray*}
   \check{N}_{\varepsilon_1}(x^*)=\{x\mid g(x)\leq0,~\|x\|_0\leq K,~\|x-x^*\|<\varepsilon_1\} \end{eqnarray*}
  and
 \begin{eqnarray*}
  \hat{N}_{\varepsilon_1}(x^*)=\{x\mid g(x)\leq0,~x_i=0, {\rm~for~}i\in  I(x^*),~
  \|x-x^*\|<\varepsilon_1\}.
 \end{eqnarray*}
  We can infer from (\ref{local_eps}) that $ \check{N}_{\varepsilon_1}(x^*)=\hat{N}_{\varepsilon_1}(x^*)$ for
  all $0<\varepsilon_1\leq \varepsilon_0$. Thus,
   $x^*$ is also a local optimal solution  of  problem $\min\{f(x)\mid g(x)\leq 0,\|x\|_0\leq S\}$. \endproof

We can infer from Theorem \ref{local_thm2} that the optimality conditions of problem  ${\rm(P}_{I(x^*)})$ are also  local
optimality conditions of problem  $\min\{f(x)\mid g(x)\leq 0,\|x\|_0\leq S\}$. And the optimality conditions of problem  ${\rm(P}_{I(x^*)})$ are easy to identify for it is a convex problem.

 Notice that testing
the feasibility of (QP) is already NP-complete when $X=\{x\mid Ax\leq b\}$ and $A$ has three rows \cite{bienstock96}.
 So it is a very tough job to find a feasible solution of the problem. Thus it is not easy to find a local optimal solution.

\begin{rem}
Constraint $x\geq 0$ can seen as  part of the constraints $g(x)\geq 0$. Then all the properties mentioned above are still true for problem (P).
\end{rem}

\section{Reformulation of the cardinality constraint}
    In this section, we propose an equivalent reformulation of problem (P).
   Let $x_{[k]}$ denote the $k$-th largest entry of $x=(x_1,\ldots,x_n)^T$. Then, when $x\geq 0$, $\|x\|_0\le K$ is equivalent to $x_{[K+1]}=0$
   or $ x_{[K+1]}\leq 0$.
Note that
\begin{eqnarray*}
    x_{[K+1]}=\sum_{i=1}^{K+1}x_{[i]}-\sum_{i=1}^{K}x_{[i]}=\psi_{K+1}(x)-\psi_{K}(x),
\end{eqnarray*}
   where $\psi_{K}(x)=\sum_{i=1}^{K}x_{[i]}$ for $K=1,\ldots,n$. The following is then obvious:
  \begin{align*}
  \psi_{K}(x)~
  =~\max_{y}~~
  &{\sum_{i=1}^{n}y_{i}x_{i}}\\
  {\rm s.t.}~~~
  &\sum_{i=1}^ny_{i}\leq K,\\
  &y_{i}\in \{0,1\},\\
  &x\geq0.
   \end{align*}

 As $\psi_{K+1}(x)-\psi_{K}(x)\leq 0$ is equivalent to  $x_{[i]}=0$ ( $i=K+1,\ldots,n$)  for all $x\geq 0$,
 $\psi_{K+1}(x)-\psi_{K}(x)\leq 0$ is then equivalent to $\sum_{i=1}^{n}x_{i}-\psi_{K}(x)\le 0$
 when $x\geq 0$. So the cardinality constraint $\|x\|_0\leq K$ is equivalent to
   \begin{eqnarray*}
    &&\sum_{i=1}^{n}x_{i}- \max_{y}~{\sum_{i=1}^{n}y_{i}x_{i}}\le 0,\\
    &&\sum_{i=1}^ny_i\leq K;\\
    &&y_{i}\in \{0,1\},\\
    &&x\geq0.
   \end{eqnarray*}
   Replacing the constraints $\sum_{i=1}^{n}x_{i}- \max_{y}~{\sum_{i=1}^{n}y_{i}x_{i}}\le 0$ and $y_{i}\in \{0,1\}$ with $\sum_{i=1}^{n}x_{i}- {\sum_{i=1}^{n}y_{i}x_{i}}\le 0$ and $0\leq y_{i}\leq 1$ gives rise to the following reformulation of problem(P):
     \begin{align*}
     {\rm(RP)}~~~~
      \min~~
      &f(x)\\
      {\rm s.t.}~~~
      &g(x)\le 0,~x\ge 0,\\
      &\sum_{i=1}^{n}x_{i}- {\sum_{i=1}^{n}y_{i}x_{i}}\le 0,\\
      &\sum_{i=1}^ny_i\leq K;\\
      &0\leq y_{i}\leq 1.
   \end{align*}
 We study in the following two subsections the global and local optimality conditions of  problem${\rm(RP)}$, respectively.

\subsection{Global optimality}
The following theorem reveals that the relaxation ${\rm(RP)}$ is essentially a reformulation of our primal problem (P) as they have the same global optimal solution and the same optimal value.
\begin{thm}\label{thm_glob}
 Problem ${\rm(RP)}$ is equivalent to  problem ${\rm (P)}$ for they have the same global optimal solutions and optimal value.
\end{thm}
\proof
We only need to prove that two conditions $x\geq 0$ and $\|x\|_0\le K$ are equivalent to
\begin{eqnarray*}
 &&x\geq 0,~~\sum_{i=1}^{n}x_{i}- {\sum_{i=1}^{n}y_{i}x_{i}}\le 0,
 \end{eqnarray*}
 where $\sum_{i=1}^ny_i\leq K$ and $0\leq y_{i}\leq 1$. Note that in the case of $x\geq0$, $\|x\|_0\le K$ is equivalent to
 \begin{eqnarray*}
 &&\sum_{i=1}^{n}x_{i}- {\sum_{i=1}^{n}y_{i}x_{i}}\le 0,\\
 &&\sum_{i=1}^ny_i\leq K,\\
 &&y_{i}\in\{0,1\}.
   \end{eqnarray*}
 Without  loss of generality, assume $y_{i}=1$ for $i=1,\ldots,l$, $0<y_{i}<1$ for $i=l+1,\ldots,m$, and $y_{i}=0$ for $i=m+1,\ldots,n$. Then, as $x\geq0$,
    \begin{eqnarray*}
   \sum_{i=1}^{n}x_{i}- {\sum_{i=1}^{n}y_{i}x_{i}}
   =\sum_{i=l+1}^{m}(1-y_{i})x_{i}+\sum_{i=m+1}^{n}x_{i}
   \leq0
    \end{eqnarray*}
   implies
   \begin{eqnarray*}
   \sum_{i=l+1}^{n}x_{i}=0,
   \end{eqnarray*}
   which, together with $\sum_{i=1}^ny_i\leq K$, yields $\|x\|_0\leq K$.\\
\hspace*{6mm}Conversely, if $\|x\|_0\leq K$,  let $y_{i}=1$ when $x_{i}>0$ and $y_{i}=0$ when $x_{i}=0$. Then it is easy to see that $\sum_{i=1}^{n}x_{i}- {\sum_{i=1}^{n}y_{i}x_{i}}\le 0$ and $\sum_{i=1}^ny_i\leq K$.
\endproof\\
\hspace*{6mm}Based on Theorem \ref{thm_glob}, we can conclude that the constraint set $\{x|\|x\|_0\leq K,~x\ge 0\}$ is equivalent to $\{x|\sum_{i=1}^{n}x_{i}- {\sum_{i=1}^{n}y_{i}x_{i}}\le 0,\sum_{i=1}^ny_i\leq K, 0\leq y_i\leq 1, i=1,\ldots,n,~x\ge0\}$. Thus, problem ${\rm RP{(x,y)}}$ is equivalent to the original  problem ${\rm (P)}$ in the sense that they have the same global solution.
 The following corollary can be obtained immediately form Theorem \ref{thm_glob}.
\begin{cor}\label{cor_golb}
 Suppose that $(\tilde{x},\tilde{y})$ is a feasible solution of problem {\rm(RP)}.
If  $\|\tilde{x}\|_0=K$, then $\tilde{y}_i\in\{0,1\}$, $\sum_{i=1}^n\tilde{y}_i=K$, and $\tilde{y}_i=1$ when $\tilde{x}_i>0$; $\tilde{y}_i=0$ when $\tilde{x}_i=0$.
 \end{cor}

\subsection{Local optimality}

It is easy to see that if $x\geq 0$, $x_i=0$ for $x_i\in I(z)$ is equivalent to $\sum_{i\in I(z)}x_i=0$. The problem  ${\rm(P}_{I(z)})$ can be restated as
\begin{align*}
     {\rm(\tilde{P}}_{I(z)})~~~~
     \min~~
     &f(x)\\
     {\rm s.t.}~~~
     &g(x)\leq0,\\
     &x\geq 0,\\
     &\hspace{-5pt}\sum_{i\in I(z)}x_i=0,
\end{align*}
We have the following theorem for local optimality.
\begin{thm}\label{local_kkt}
 Suppose ($x^*,y^*$) is a  KKT point of ${\rm(RP)}$ with $\|x^*\|_0=K$.  Then $x^*$  also satisfies the KKT conditions of problem ${\rm(\tilde{P}_{I({x^*})})}$, and  $x^*$ is a local optimal solution of problem $\rm{(P)}$.
\end{thm}
\proof
With  $\|x^*\|_0=K$, we have $y_i^*\in\{0,1\}$ and $\sum_{i=1}^ny^*_i=K$ from Corollary \ref{cor_golb}.
If ($x^*,y^*$) is a  KKT point of ${\rm(RP)}$,  then
\begin{eqnarray}
&&\nabla f({x}^*)+\sum_{i=1}^m v_i \nabla g_i({x}^*)+\delta(\sum_{i=1}^n(1-y_i^*)e_i)+\sum_{i=1}^n\tau_ie_i=0,\label{kkt-3-1}\\
&&v_i g_i({x^*})=0,\label{kkt-3-2}\\
&&\tau_i{x^*}_i=0,~i=1,\ldots,n.\label{kkt-3-3}\\
&&\delta(\sum_{i=1}^nx^*_i-\sum_{i=1}^ny^*_ix^*_i)=0\label{kkt-3-4}\\
&&-\delta(\sum_{i=1}^nx_i^*e_i)+\rho e+\sum_{i=1}^n \varrho_ie_i-\sum_{i=1}^n\sigma_ie_i=0, \\
&&\rho(\sum_{i=1}^ny^*_i-K)=0,\label{kkt-3-5}\\
&& \varrho_i y^*_i=0,~~i=1,\ldots,n,\\
&&\sigma_i(y^*_i-1)=0,~~i=1,\ldots,n, \end{eqnarray}
where  $ v_i \in\Re^1_+$, $\tau_i\in\Re^1_+$, $\delta\in\Re^1_+$, $\rho\in\Re^1_+$, $\varrho_i\in\Re^1_+$,
$\sigma_i\in\Re^1_+$, $e$ denotes the all-one vector and $e_i$ denotes the $i$-th coordinate vector.
Since $y_i^*\in\{0,1\}$, $\sum_{i=1}^ny^*_i=K$  and  $\|x^*\|_0=K$,  we have
 \begin{eqnarray}
 \sum_{i=1}^nx^*_i-\sum_{i=1}^ny^*_ix^*_i=\sum_{i\in I(x^*)}x_i^*=0 \label{kkt-4-1}
 \end{eqnarray}
and
 \begin{eqnarray}
\sum_{i=1}^n(1-y_i^*)e_i=\sum_{i\in I(x^*)}e_i.\label{kkt-5-1}
 \end{eqnarray}
Combining  (\ref{kkt-3-1}), (\ref{kkt-3-2}),(\ref{kkt-3-3}), (\ref{kkt-4-1}) and (\ref{kkt-5-1}) yields
 \begin{eqnarray}
 &&\nabla f({x}^*)+\sum_{i=1}^m v_i \nabla g_i({x}^*)+\delta(\sum_{i\in I(x^*)}e_i)+\sum_{i=1}^n\tau_ie_i=0,\label{kkt-2-1}\\
 &&v_i g_i({x^*})=0,\label{kkt-2-2}\\
 &&\tau_i{x^*}_i=0,~i=1,\ldots,n,\label{kkt-2-3}\\
 &&\delta(\sum_{i\in I(x^*)}x_i^*)=0\label{kkt-2-4},
 \end{eqnarray}which means  that $x^*$  also satisfies the KKT conditions of problem ${\rm(\tilde{P}}_{I({x^*})})$, i.e. $x^*$ is global optimal solution of problem ${\rm(\tilde{P}}_{I({x^*})})$.
  From Theorem
 \ref{local_thm2},  $x^*$ is also a local optimal solution of Problem (P).
\endproof

\bigskip
Theorem \ref{local_kkt} means that when we get a KKT point of  ${\rm(RP)}$ by some algorithm, we essentially get a local optimal
solution of problem (P) under certain conditions.
\begin{rem}
 Let $z_i=1-y_i$. The transformation  ${\rm(RP)}$  can be equivalently transformed into an optimization problem with complementarity constraints, which is  discussed  in \cite{Burdakov2014} and \cite{Feng2013}.
Some properties about the ${\rm(RP)}$ problem can also be found in  \cite{Burdakov2014} by taking this transformation.
By the transformation, Theorem 3 can be transformed
into Theorem 3.2 in  \cite{Burdakov2014}. Corollary 2 can be transformed into Proposition 3.5 in  \cite{Burdakov2014}.
The contribution of our method is that
 construct an algorithm
based on the model in the next section.  This algorithm can  provide a local optimal solution of problem (P), and this local optimal solution  exactly satisfies the cardinality constraint as we shown in Theorem 1.
\end{rem}

\section{A successive convex approximation method}
 In this section, we first  propose a successive convex approximation method for the reformulation problem ${\rm(RP)}$ by constructing a sequence of convex subproblems. We then establish the convergence of the method to a KKT point of ${\rm(RP)}$.
 Based on Theorem\ref{local_thm2}  and \ref{local_kkt}, we get a local optimal solution of problem (P)  under certain conditions.
\subsection{Approximation method}
 Let $(\bar{x},\bar{y})$ be a feasible solution to problem ${\rm(RP)}$.
  We use the first order  Taylor expansion to approximate ${\sum_{i=1}^{n}y_{i}x_{i}}$:
   \[
   {\sum_{i=1}^{n}y_{i}x_{i}}\approx \sum_{i=1}^{n}\bar{y_{i}}\bar{x_{i}}+\bar{\xi}^T(x-\bar{x},y-\bar{y}),~~~~\forall x,y\in \mathcal{R}^n,
   \]
   where $\bar{\xi}=\nabla f(\bar{x},\bar{y})=(\bar{y}^T,\bar{x}^T)^T$. Then
 \begin{eqnarray*}
   \sum_{i=1}^{n}x_{i}-{\sum_{i=1}^{n}y_{i}x_{i}}&\approx& \sum_{i=1}^{n}x_{i}-( \sum_{i=1}^{n}\bar{y_{i}}\bar{x_{i}}+\bar{\xi}^T(x-\bar{x},y-\bar{y}))\\
   &=&\sum_{i=1}^{n}x_{i}-( \sum_{i=1}^{n}\bar{y_{i}}\bar{x_{i}}+(\bar{y}^T(x-\bar{x})+\bar{x}^T(y-\bar{y}))\\
   &=&(e-\bar{y})^Tx-\bar{x}^Ty+\sum_{i=1}^{n}\bar{x_i}\bar{y_i}.
 \end{eqnarray*}
So a convex approximation model of ${\rm(RP)}$ at $(\bar{x},\bar{y})$   can be presented as follows:
      \begin{align}
     {\rm(AP}(\bar{x},\bar{y}))~~~~
     \min~~
     &f(x)\\
     {\rm s.t.}~~~
     &g(x)\leq 0,\\
     &\phi_{(\bar{x},\bar{y})}(x,y)\triangleq(e-\bar{y})^Tx-\bar{x}^Ty+\sum_{i=1}^{n}\bar{x}_i\bar{y}_i\leq 0,\label{apm1}\\
     &\sum_{i=1}^ny_i\leq K,\\
     &0\leq y_{i}\leq 1,~i=1,\ldots,n,\\
     &x\geq 0.
   \end{align}
In the following lemma, we  present the relationship between optimal solutions of (${\rm AP}(\bar{x},\bar{y})$) and local
optimal solutions of problem (P) under certain conditions.

\begin{lem}\label{lem_opt}
If  $(\bar{x},\bar{y})$ is an optimal solution to problem
(${\rm AP}(\bar{x},\bar{y})$) with $\|\bar{x}\|_0=K$,
then $\bar{x}$ is  a local optimal solution of problem  ${\rm(P})$.
\end{lem}
\proof Because $(\bar{x},\bar{y})$ is feasible to (${\rm AP}(\bar{x},\bar{y}))$), we have
\[
 (e-\bar{y})^T\bar{x}-\bar{x}^T\bar{y}+\sum_{i=1}^{n}\bar{x}_i\bar{y}_i=(e-\bar{y})^T\bar{x}=0.
\]
Because of $\|\bar{x}\|_0=K$, we have $\bar{y}_i\in\{0,1\}$ for $i=1,\ldots,n$,
and $\bar{y}_i=1$ when $\bar{x}_i>0$; $\bar{y}_i=0$ when $\bar{x}_i=0$.
Thus $(e-\bar{y})^Tx=\sum_{i\in I(\bar{x})}x_i$
and $\bar{x}$ is feasible to problem ${\rm(\tilde{P}}_{I(\bar{x})})$.
Note that if $x$ is feasible to problem ${\rm(\tilde{P}}_{I(\bar{x})})$, then
$(x,\bar{y})$ is feasible to problem ${\rm(AP}(\bar{x},\bar{y}))$.
Thus
$\bar{x}$ is also an optimal solution of problem ${\rm(\tilde{P}}_{I(\bar{x})}) $. According  to Theorem \ref{local_thm2}, $\bar{x}$ is a local optimal solution of problem (P).\endproof

\subsection{SCA method}
In this section, we will develop a successive convex approximation (SCA) method for the problem we study. In this method, when we treat
${\rm(AP}(\bar{x},\bar{y}))$ as a subproblem, we  need a feasible solution of problem (P) to start with.  Notice  that testing
the feasibility of (QP) is already NP-complete when $X=\{x\mid Ax\leq b\}$ and $A$ has three rows \cite{bienstock96}.
 So it is also a very tough job to find a feasible solution of problem (P). Here  we use  the penalty method
 to handle  constraint (\ref{apm1}).  By introducing a  penalty parameter  $\mu$,  we consider the following convex subproblem:
 \begin{align*}
     {\rm(AP_{\mu}(\bar{x},\bar{y}))}~~~~
     \min~~
     &f(x)+\mu\left[(e-\bar{y})^Tx-\bar{x}^Ty+\sum_{i=1}^{n}\bar{x}_i\bar{y}_i\right]\\
     {\rm s.t.}~~~
     &g(x)\leq 0,\\
     &\sum_{i=1}^ny_i\leq K,\\
     &0\leq y_{i}\leq 1,~i=1,\ldots,n,\\
     &x\geq 0.
   \end{align*}

A successive approximation method for ${\rm(RP)}$ is presented in Algorithm \ref{algo1}.
An initial solution $x^0$ in Step 1 of Algorithm \ref{algo1} can be
obtained by solving the following convex quadratic programming:
\begin{align*}
  \min~~
  &f(x)\\
  {\rm s.t.}~~~
  &g(x)\leq 0,\\
  &x\geq0.
\end{align*}
Our theoretical and computational results show effectiveness of this proposed scheme.

\begin{algo}[SCA method for ${\rm(RP)}$]
\label{algo1}
\vspace{5pt}
\begin{enumerate}[label=\textit{Step \arabic*}:, leftmargin=*]
\item Choose an initial  positive value for parameter $\mu=\mu_0$, a positive value for $\varrho$, and a small positive value for stopping parameter $\epsilon$.
Select   $x^0\in X=\{x\mid g(x)\leq0,~x\ge0\}$.  If $\|x^0\|_0\leq K$, stop; else
set $k=0$ and
\begin{align*}
   y^0_i =  \begin{cases}
                    1 &{\rm~if~}x^0_i\ge x^0_{[K]},  \\
                    0 &{\rm~otherwise},  \\
              \end{cases}~~~~~~
              i=1,\ldots,n.
\end{align*}

\item Solve the convex subproblem ${\rm AP}_\mu{({x}^k,{y}^k)}$ to get $(x^{k+1},y^{k+1})$.

\item If $\|x^{k+1}-x^k,y^{k+1}-y^k\|\le \epsilon$ and $\|x^{k+1}\|\leq K$, stop.

\item If $(e-{y^{k}})^Tx^{k+1}>0$,
set
\begin{align*}
   y^{k+1}_i =  \begin{cases}
                    1 &{\rm~if~}x^{k+1}_i\ge x^{k+1}_{[K]},  \\
                    0 &{\rm~otherwise},  \\
              \end{cases}~~~~~~
              i=1,\ldots,n.
\end{align*}

\item If $\|y^{k+1}\|_0>K$, set
\begin{align*}
   y^{k+1}_i =  \begin{cases}
                    1 &{\rm~if~}x^{k}_i\ge x^{k}_{[K]},  \\
                    0 &{\rm~otherwise},  \\
              \end{cases}~~~~~~
              i=1,\ldots,n.
\end{align*}
%

\item Set $k:=k+1$, $\mu=\varrho \mu$ and go to Step 2.

\end{enumerate}
\end{algo}

%
%
%
%
%
%

We assume $\|{x}'\|_0>K$ for all ${x}'\in \arg \min\{f(x)\mid g(x)\leq 0,~x\geq 0\}$.
We also assume that the inner of set $\{x|g(x)\leq 0,x\geq 0,~\|x_0\|\leq K\}$ is nonempty.
The assumptions are not so difficult to satisfy. For example, the cardinality constraint portfolio problems we  use as the computational cases in the next section satisfy these assumptions.
\begin{rem}
Based on our assumption mentioned above and Theorem 1, the vector $x^{k+1}$ could not be too sparse, i.e., $\|x^{k+1}\|<K$.
In Step 4 and Step 5, if the $K$-th largest entry of $x^{k+1}$ or $x^k$
is not unique, the elements of   $y^{k+1}$ corresponding the $k$ largest elements of $x^{k+1}$ or $x^k$ would be 1, and the other elements of $y^{k+1}$ would be zero.
\end{rem}

\begin{lem}\label{lem_algo}
 Let $(x^k,y^k)$ ($k=1,2,\ldots$) be the sequence generated by  Algorithm 1. Then $y^k_i\in\{0,1\}$, $\sum_{i=1}^ny^k_i=K$, $\sum_{i=1}^{n}{x_i^{k-1}}{y_i^{k-1}}-{(x^{k-1})}^Ty^k=0$,
 $(e-{y^{k-1}})^Tx^k-{(x^{k-1})}^Ty^k+\sum_{i=1}^{n}{x_i^{k-1}}{y_i^{k-1}}\geq 0$, $\|x^k\|_0\geq K$,
and
 \begin{align*}
    y^{k}_i =  \begin{cases}
                     1 &{\rm~if~}x^{k}_i\ge x^{k}_{[K]},  \\
                     0 &{\rm~otherwise}.  \\
                \end{cases}~~~~~~
           i=1,\ldots,n.
 \end{align*}
\end{lem}
 \proof
 We  prove these properties by induction on $k$. Firstly,
 when $k=1$, the objective function of the  convex subproblem ${\rm(AP_\mu({x^{0}},{y^{0}}))}$ is
 \begin{eqnarray*}
  F(x,y)=f(x)+\mu[(e-{y^{0}})^Tx-{(x^{0})}^Ty+\sum_{i=1}^{n}{x_i^{0}}{y_i^{0}}].
 \end{eqnarray*}
 It is obvious that ${(x^{0})}^Ty$ should be as big as possible when  minimizing the objective function.
 Suppose the optimal solution of ${\rm(AP_\mu({x^{0}},{y^{0}}))}$ is $(x^1,y^1)$.
Then we must have
  ${(x^{0})}^Ty^{1}={(x^{0})}^Ty^{0}=\sum_{i=1}^Kx_{[i]}^{0}$ and $\sum_{i=1}^ny_i^{1}=K$
because $\sum_{i=1}^nx_{i}^{0}y_i^0=\sum_{i=1}^Kx_{[i]}^{0}$, $x\geq 0$, $0\leq y\leq e$ and  $\sum_{i=1}^ny_i\leq K$.
Thus we have $\|y^1\|_0\geq K$ because  $0\leq y^1\leq e$ and  $\sum_{i=1}^ny^1_i= K$.
 Then the objective function of problem ${\rm(AP_\mu({x^{0}},{y^{0}}))}$ is equivalent to  $F(x,y)=f(x)+\mu[(e-{y^{0}})^Tx]$.
It is obvious that $(e-{y^{0}})^Tx\geq 0$ because $x\geq 0$ and $0\leq y^0\leq e$.   If  $(e-{y^{0}})^Tx^1=0$,  we have $\|x^1\|_0\leq K$ because $0\leq y^0_i\leq 1$ and $\sum_{i=1}^ny^0_i=K$. According to
 Corollary \ref{cor_th1}, we can get that
$\|x^1\|_0=K$.
If  $(e-{y^{0}})^Tx^1>0$, we will prove that  $\|x^1\|_0\geq K$. By contradiction, suppose $\|x^1\|_0<K$. Let $\hat{x}^1$ be the optimal solution of the following problem,
\begin{align*}
    ~~~~
    \min~~
    &f(x)\\
    {\rm s.t.}~~~
    &g(x)\leq 0,\\
    &x\geq 0,\\
    &(e-{y^{0}})^Tx=\sum_{i\in I(y^0)}x_i=0.
  \end{align*}
 Then  $\|\hat{x}^1\|_0=K$ based on  Corollary \ref{cor_th1}.
  So, $\hat{x}^1$ is a local optimal solution of problem (P) according to Theorem \ref{local_thm2}.
   Since $\hat{x}^1$ is a feasible solution of problem (P),  $(\hat{x}^1,y^1)$  is feasible to  ${\rm(AP_\mu({x^{0}},{y^{0}}))}$.
 Then
\begin{eqnarray*}
&&f(x^1)+\mu[(e-{y^{0}})^Tx^1-{(x^{0})}^Ty^1+\sum_{i=1}^{n}{x_i^{0}}{y_i^{0}}]=f(x^1)+\mu[(e-{y^{0}})^Tx^1]\\
&&\leq f(\hat{x}^1)+\mu[(e-{y^{0}})^T\hat{x}^1]= f(\hat{x}^1).
\end{eqnarray*}
Thus, $f(x^1)\leq f(\hat{x}^1)$, which contradicts Corollary \ref{cor_th1}. Therefore,  $\|x^1\|_0\geq K$.

Now we  prove that $y_i^1\in\{0,1\}$ ($i=1,\ldots,n$) in Algorithm 1. It is obvious that if $\|y^1\|_0=K$, then $y_i^1\in\{0,1\}$ because $\sum_{i=1}^ny_i^1=K$ and $0\leq y^1\|_0\leq e$.
In this case, if $(e-{y^{0}})^Tx^1=0$, we have $\|x^1\|_0=K$. Then we can infer that $x^1_i>0$
when $y^0_i=1$  and  $x^1_i=0$ when $y^0_i=0$ from $y_i^1\in\{0,1\}$, $\sum_{i=1}^ny_i^1=K$, $\|y^1\|_0=K$ and  $\|x^1\|_0=K$. In the case of $\|y^1\|_0=K$, i.e. $y_i^1\in\{0,1\}$, if $(e-{y^{0}})^Tx^1>0$, then we have
\begin{align*}
   y^{1}_i =  \begin{cases}
                    1 &{\rm~if~}x^{1}_i\ge x^{1}_{[K]},  \\
                    0 &{\rm~otherwise},  \\
              \end{cases}~~~~~~
              i=1,\ldots,n,
\end{align*}
 from Step 4 of Algorithm 1.

 If  $\|y^1\|_0>K$, then we set
\begin{align*}
   y^{1}_i =  \begin{cases}
                    1 &{\rm~if~}x^{0}_i\ge x^{0}_{[K]},  \\
                    0 &{\rm~otherwise},  \\
              \end{cases}~~~~~~
         i=1,\ldots,n,
\end{align*}
 in Step 5 of Algorithm 1. Thus,  ${(x^{0})}^Ty^{1}=\sum_{i=1}^nx_{i}^{0}y_i^0$ and $y_i^1\in\{0,1\}$.
 Therefor we have $y^1_i\in\{0,1\}$, $\sum_{i=1}^ny^1_i=K$, $\sum_{i=1}^{n}{x_i^{0}}{y_i^{0}}-{(x^{0})}^Ty^1=0$,
 $(e-{y^{0}})^Tx^1-{(x^{0})}^Ty^1+\sum_{i=1}^{n}{x_i^{0}}{y_i^{0}}\geq 0$, $\|x^1\|_0\geq K$,
 and
 \begin{align*}
    y^{1}_i =  \begin{cases}
                     1 &{\rm~if~}x^{1}_i\ge x^{1}_{[K]},  \\
                     0 &{\rm~otherwise}.  \\
                \end{cases}~~~~~~
           i=1,\ldots,n.
 \end{align*}                       \\

Suppose that these properties  hold before iteration $k-1$. Now we consider the convex subproblem ${\rm(AP_\mu({x^{k-1}},{y^{k-1}}))}$.
The objective function of the  convex subproblem ${\rm(AP_\mu({x^{k-1}},{y^{k-1}}))}$ is
\begin{eqnarray*}
 F(x,y)=f(x)+\mu[(e-{y^{k-1}})^Tx-{(x^{k-1})}^Ty+\sum_{i=1}^{n}{x_i^{k-1}}{y_i^{k-1}}].
\end{eqnarray*}
It is obvious that ${(x^{k-1})}^Ty$ should be as big as possible when  minimizing the objective function.
Suppose the optimal solution of ${\rm(AP_\mu({x^{k-1}},{y^{k-1}}))}$ is $(x^k,y^k)$. Then
we must have ${(x^{k-1})}^Ty^{k}=\sum_{i=1}^Kx_{[i]}^{k-1}$ and $\sum_{i=1}^ny_i^{k}=K$
because $\sum_{i=1}^nx_{i}^{k-1}y_i^{k-1}=\sum_{i=1}^Kx_{[i]}^{k-1}$, $x\geq 0$, $0\leq y\leq e$ and  $\sum_{i=1}^ny_i\leq K$.
 Then we have $\|y^k\|_0\geq K$ and the objective function of problem ${\rm(AP_\mu({x^{k-1}},{y^{k-1}}))}$ is equivalent to  $F(x,y)=f(x)+\mu[(e-{y^{k-1}})^Tx]$.
It is obvious that $(e-{y^{k-1}})^Tx\geq 0$ because $x\geq 0$ and $0\leq y\leq e$.
If  $(e-{y^{k-1}})^Tx^k=0$,  we have $\|x^k\|_0\leq K$ for $0\leq y^{k-1}_i\leq 1$ and $\sum_{i=1}^ny^{k-1}_i=K$. According to Corollary \ref{cor_th1},
$\|x^k\|_0=K$ must hold.

If  $(e-{y^{k-1}})^Tx^k>0$, we  prove  $\|x^k\|_0\geq K$. By contradiction, suppose $\|x^k\|_0<K$. Let
$\hat{x}^k$ be the optimal solution of the following  problem,
 \begin{align*}
     ~~~~
     \min~~
     &f(x)\\
     {\rm s.t.}~~~
     &g(x)\leq 0,\\
     &x\geq 0,\\
     &(e-{y^{k-1}})^Tx=\sum_{i\in I(y^{k-1})}x_i=0.
\end{align*}
Then  $\|\hat{x}^k\|_0=K$ due to  Corollary \ref{cor_th1}.
So, $\hat{x}^k$ is a local optimal solution of problem (P) according to Theorem \ref{local_thm2}.
Since $\hat{x}^k$ is a feasible solution of problem (P), $(\hat{x}_k,y^k)$  is feasible for
${\rm(AP_\mu({x^{k-1}},{y^{k-1}}))}$.
 Then
\begin{eqnarray*}
&&f(x^k)+\mu[(e-{y^{k-1}})^Tx^k-{(x^{k-1})}^Ty^k+\sum_{i=1}^{n}{x_i^{k-1}}{y_i^{k-1}}]=f(x^k)+\mu[(e-{y^{k-1}})^Tx^k]\\
&&\leq f(\hat{x}^k)+\mu[(e-{y^{k-1})^T\hat{x}^k}]= f(\hat{x}^k).
\end{eqnarray*}
Thus, $f(x^k)< f(\hat{x}^k)$ which contradicts Corollary \ref{cor_th1}. So, $\|x^k\|_0\geq K$.

Now we  prove $y_i^k\in\{0,1\}$ ($i=1,\ldots,n$) in Algorithm 1.
It is obvious that if $\|y^k\|_0=K$, then $y_i^k\in\{0,1\}$ because $\sum_{i=1}^ny_i^k=K$ and $0\leq y^1\leq e$. In this case, if $(e-{y^{k-1}})^Tx^k=0$, we have $\|x^k\|_0=K$. Then we can infer that $x^k_i>0$ when $y^{k-1}_i=1$  and  $x^k_i=0$ when $y^{k-1}_i=0$ from
 $y^{k-1}_i\in\{0,1\}$, $\sum_{i=1}^ny^{k-1}_i=K$, $\|x^k\|_0=K$ and $(e-{y^{k-1}})^Tx^k=0$.  In the case of $\|y^k\|_0=K$, i.e. $y_i^k\in\{0,1\}$,
 if $(e-{y^{k-1}})^Tx^k>0$, from Step 4 of Algorithm 1, we have
\begin{align*}
   y^{k+1}_i =  \begin{cases}
                    1 &{\rm~if~}x^{k+1}_i\ge x^{k+1}_{[K]},  \\
                    0 &{\rm~otherwise},  \\
              \end{cases}~~~~~~
              i=1,\ldots,n.
\end{align*}

  On the other side, if
  $\|y^k\|_0>K$,   then we set
 \begin{align*}
    y^{k}_i =  \begin{cases}
                     1 &{\rm~if~}x^{k-1}_i\ge x^{0}_{[K]},  \\
                     0 &{\rm~otherwise},  \\
               \end{cases}~~~~~~
          i=1,\ldots,n.
 \end{align*}
  in Step 5 of Algorithm 1. Thus ${(x^{k-1})}^Ty^{k}=\sum_{i=1}^nx_{i}^{k-1}y_i^{k-1}$ and $y_i^{k}\in\{0,1\}$.
 Therefore we have $y^k_i\in\{0,1\}$, $\sum_{i=1}^ny^k_i=K$, $\sum_{i=1}^{n}{x_i^{k-1}}{y_i^{k-1}}-{(x^{k-1})}^Ty^k=0$,
  $(e-{y^{k-1}})^Tx^k-{(x^{k-1})}^Ty^k+\sum_{i=1}^{n}{x_i^{k-1}}{y_i^{k-1}}\geq 0$, $\|x^k\|_0\geq K$,
 and
  \begin{align*}
     y^{k}_i =  \begin{cases}
                      1 &{\rm~if~}x^{k}_i\ge x^{k}_{[K]},  \\                                                   0 &{\rm~otherwise}.  \\
                 \end{cases}~~~~~~
            i=1,\ldots,n.
  \end{align*}

\endproof

We now establish the convergence of Algorithm 1 to  a local optimal solution of problem (P).

\begin{thm}\label{lemm_3}Let $\epsilon=0$.
{\rm(i)} If the algorithm stops at Step 3 in the iteration, then $x^{k}$ is a local optimal solution of problem (P).

{\rm(ii) }If the algorithm generates
an infinite sequence $(x^k, y^k)$, then any accumulation point $(x^*, y^*)$ is a global optimal solution
of problem ${\rm (P_{I(x^*)})}$, i.e. $x^*$ is a local optimal solution of Problem (P).

\end{thm}
\proof
(i) If the algorithm stops at Step 3, then  ($x^k,y^k$) solves subproblem
${\rm(AP_\mu({x^k},{y^k}))}$, and
\[
(e-{y^k})^Tx^k-{x^k}^T{y^k}+\sum_{i=1}^{n}x^k_iy^k_i=(e-{y^k})^Tx^k=0
\]
 because $\|x^k\|_0=K$
and
\begin{align*}
   y^{k}_i =  \begin{cases}
                    1 &{\rm~if~}x^{k}_i\ge x^{k}_{[K]},  \\
                    0 &{\rm~otherwise},  \\
               \end{cases}~~~~~~
          i=1,\ldots,n,
\end{align*}
from Lemma \ref{lem_algo}.
So $(x^{k},y^k)$ is also the optimal solution of problem ${\rm AP({x^k},{y^k})}$.  Thus from Lemma \ref{lem_opt}, $x^{k}$ is
a local optimal solution of problem (P).

\hspace*{6mm}(ii) Suppose $(x^k,y^k)$ is the optimal solution of ${\rm(AP_{\mu}({x^{k-1}},{y^{k-1}}))}$. Then  $y^k_i\in\{0,1\}$ and $\sum_{i=1}^ny^k_i=K$ according to Lemma \ref{lem_algo}.
Let $\tilde{x}$ be the optimal solution
 of  the following convex problem:
 \begin{align*}
     ~~~~
     \min~~
     &f(x)\\
     {\rm s.t.}~~~
     &g(x)\leq 0,\\
     &x\geq 0,\\
     &(e-{y^{k}})^Tx=\sum_{i\in I(y^k)}x_i=0.
\end{align*}
Then $I(\tilde{x})=I(y^{k})$ and $\|\tilde{x}\|_0=K$  based on Corollary \ref{cor_th1}. So $\tilde{x}$ is a local
optimal solution of problem (P) from Theorem \ref{local_thm2}. It is obvious that $(\tilde{x}, y^{k})$ is feasible for  ${\rm(AP_{\mu}({x^k},{y^k}))}$, then
\begin{eqnarray}
&&F(x^{k+1},y^{k+1})=f(x^{k+1})+\mu_k[(e-{y^{k}})^Tx^{k+1}-{x^k}^Ty^{k+1}+\sum_{i=1}^{n}x_i^ky_i^k]\\
&=&f(x^{k+1})+\mu_k[(e-{y^{k}})^Tx^{k+1}]\leq f(\tilde{x})+\mu_k[(e-{y^{k}})^T\tilde{x}]=f(\tilde{x}),~k=1,2,\ldots. \label{local_algo}
\end{eqnarray}
Hence,
\begin{eqnarray}
0\leq (e-{y^{k}})^Tx^{k+1}\leq \frac{f(\tilde{x})-f(x^{k+1})}{\mu_k},
\end{eqnarray}
because of $-{x^k}^Ty^{k+1}+\sum_{i=1}^{n}x_i^ky_i^k=0$. Consider two converging subsequences $\{x^{k},y^k\}$ and $\{x^{k+1},y^{k+1}\}$. Let $(\hat{x^*},y^*)$ and $({x^*},\hat{y^*})$ be  limits of $\{x^{k},y^k\}$ and $\{x^{k+1},y^{k+1}\}$ respectively. Passing to the limit with $k\rightarrow \infty$ in
the above inequality, we conclude that $(e-{y^{k}})^Tx^{k+1}\rightarrow 0$. Then we can get that $(e-{y^{*}})^T{x^{*}}=0$,
$I(\tilde{x})=I(y^{*})$, $0\leq y^{*}\leq e$, $\sum_{i=1}^ny_i^*=K$, and $x^*\geq 0$. Thus  $\|x^*\|=K$,  $y^*_i=1$ when $x^*_i>0$ and $y^*_i=0$ when $x^*_i=0$.

Furthermore, inequality (\ref{local_algo}) implies
\begin{eqnarray}
f(x^{*})+\mu_k[(e-{y^{*}})^Tx^{*}]=f(x^{*})\leq f(\tilde{x}).
\end{eqnarray}
Together with $I(\tilde{x})=I(y^{*})$,  $x^*$ is a local optimal solution of problem (P).
\endproof \\

Now we show that the penalty parameter $\mu$ should not be too big under certain conditions.
Let $F_0=\{(x,y)|g(x)\leq 0,\sum_{i=1}^ny_i\leq K,x\geq 0, 0\leq y_{i}\leq 1,i=1,\ldots,n.\}$. Then the convex subproblem
${\rm(AP}(\bar{x},\bar{y}))$ can be restated as
       \begin{align*}
      {\rm(AP'}(\bar{x},\bar{y}))~~~~
      \min~~
      &f(x)\\
      {\rm s.t.}~~~
      &(e-\bar{y})^Tx-\bar{x}^Ty+\sum_{i=1}^{n}\bar{x}_i\bar{y}_i\le 0,\\
      &(x,y)\in F_0,
      \end{align*}
and the convex  subproblem  ${\rm(AP_{\mu}(\bar{x},\bar{y}))}$ can be restated as
\begin{align*}
    {\rm(AP'_{\mu}(\bar{x},\bar{y}))}~~~~
    \min~~
    &f(x)+\mu[(e-\bar{y})^Tx-\bar{x}^Ty+\sum_{i=1}^{n}\bar{x}_i\bar{y}_i]\\
    {\rm s.t.}~~~
    &(x,y)\in F_0.
\end{align*}
\begin{thm}\label{penal_p}
If $(\bar{x},\bar{y})$ is a global minimum of problem ${\rm(AP'}(\bar{x},\bar{y}))$ which satisfies the second order sufficient conditions of optimality with multipliers $\hat{\lambda}$ (see Theorem 3.47 in \cite{Ruszczynsk2006}). Then for every $\mu>\hat{\lambda}$, $(\bar{x},\bar{y})$ is a global
minimum of problem ${\rm(AP'_{\mu}(\bar{x},\bar{y}))}$.
\end{thm}
\proof
 Based on Theorem 6.9 of \cite{Ruszczynsk2006}, we can conclude that  $(\bar{x},\bar{y})$ is a global minimum of the following convex problem:
 \begin{align*}
     {\rm(AP''_{\mu}(\bar{x},\bar{y}))}~~~~
     \min~~
     &f(x)+\mu\max\{0,((e-\bar{y})^Tx-\bar{x}^Ty+\sum_{i=1}^{n}\bar{x}_i\bar{y}_i\}\\
     {\rm s.t.}~~~
     &(x,y)\in F_0.
 \end{align*}
It is easy to see that $(e-\bar{y})^T\bar{x}-(\bar{x})^T\bar{y}+\sum_{i=1}^{n}\bar{x}_i\bar{y}_i=e^T\bar{x}-\bar{x}^T\bar{y}\geq 0$ for $\bar{x}\geq0$ and $0\leq \bar{y}\leq e$.
 Then  $(\bar{x},\bar{y})$ is also a global minimum of problem ${\rm(AP'_{\mu}(\bar{x},\bar{y}))}$.
\endproof \\

 The existence of Lagrange multipliers and the ability of the penalty model to get exact solutions are closed related.  Based on Theorem 6.10 of \cite{Ruszczynsk2006}, we can get the following theorem easily.
\begin{thm}\label{penal_F}
Suppose that $F_0$ is a closed polyhedron and that $(\bar{x},\bar{y})$ is a global minimum of the convex problem
  ${\rm(AP''_{\mu}(\bar{x},\bar{y}))}$. If  $(\bar{x},\bar{y})$ is feasible for problem ${\rm(AP'}(\bar{x},\bar{y}))$, then there exists Lagrange multiplier $\lambda$ such that $(\bar{x},\bar{y})$ satisfies the first order necessary conditions of optimality
  for convex subproblem ${\rm(AP'}(\bar{x},\bar{y}))$. Furthermore,  $(\bar{x},\bar{y})$ is the global optimal
  solution of the convex subproblem ${\rm(AP'}(\bar{x},\bar{y}))$.
\end{thm}

\begin{rem}
Based on Theorem \ref{lemm_3}, we can get a local optimal solution which exactly satisfied the cardinality constraint under certain conditions.
 Based on Theorems \ref{penal_p} and \ref{penal_F}, we can see that the penalty parameter $\mu$ should not be  too large in Algorithm 1 under certain conditions.
If $(\bar{x},\bar{y})$ is the  global optimal solution of the convex subproblem
 ${\rm(AP'}(\bar{x},\bar{y}))$, then $e^T\bar{x}-\bar{x}^T\bar{y}=0$ and $\|\bar{x}\|_0\leq K$. If $\|\bar{x}\|_0= K$, then $\bar{x}$ is also a local optimal solution of problem $\rm{(P)}$ under certain conditions based on
Lemma \ref{lem_opt} and Theorem \ref{local_kkt}. If Algorithm 1 stops at Step 3,  then
$(x^{k},y^{k})$ is optimal for ${\rm(AP'_{\mu}({x^k},{y^k}))}$. Thus, $e^T{x^k}-{x^k}^T{y^k}=0$ and $\|{x^k}\|_0\leq K$.
If the algorithm generates
an infinite sequence $(x^k, y^k)$, then any accumulation point $(x^*, y^*)$ is optimal solution of ${\rm(AP'_{\mu}({x^*},{y^*}))}$.
Thus, $e^T{x^*}-{x^*}^T{y^*}=0$ and $\|{x^*}\|_0\leq K$.  If $\|{x^*}\|_0=K$ we can get a local optimal solution of problem $\rm{(P)}$.
Theorems \ref{penal_p} and  \ref{penal_F} show that we do not need a too large  penalty parameter $\mu$ to get a local optimal solution of ${\rm(P_1)}$ under certain conditions.
Actually we will show these results in our computational experiments in the next section.
\end{rem}

\section{Computational results}

We use limited diversified
mean-variance portfolio selection problems (see \cite{bienstock96,bonami09}) as the test problems in our computational experiments.
 The variables of the problems are all confined to be nonnegative and the continuous relaxation of the feasible set of the problem is a closed polyhedral. In our computational experiments,
we compare Algorithm 1 with two  successive  approximation (SCA) methods using ``$\ell_p$'' approximation  and exponential approximation  in \cite{chen2010a,mangasarian1999}, respectively.
The  two SCA methods are called
   ``SCA-$\ell_p$'' and ``SCA-exp'' methods,  where the convex approximation subproblems are obtained from ``$\ell_p$'' function and
exponential  function respectively in \cite{chen2010a,mangasarian1999}. The ``$\ell_p$'' function is
\[
\ell_p(x)=\sum_{i=1}^nx_i^p,
\]
where $p$ is a scalar parameter with $0<p<1$. The exponential approximation function is
\[
{\rm exp}_p(x)=\sum_{i=1}^n (1-e^{-\frac{1}{p}x_i}),
\]
where $p>0$ is a scalar parameter. The structure of the two SCA methods  using ``$\ell_p$'' approximation  or exponential approximation is shown in Algorithm 2.
\begin{algo}[SCA Method ]
\label{algo2}
\vspace{5pt}

{\rm
\begin{enumerate}[label=\textit{Step \arabic*}:, leftmargin=*]

\item Choose a positive value for parameter $\mu$ and a small positive value for stopping parameter $\epsilon$.
Select  $x^0$ satisfying $x^0\in X=\{x|g(x)\leq0\}$. Set $k=0$.

\item Solve a convex approximation subproblem obtained from  parameter $\mu$ and ``$\ell_p$'' approximation function or
 exponential approximation function to get $x^{k+1}$.

\item If $\|x^{k+1}-x^k\|\le \epsilon$, stop.

\item Set $k:=k+1$ and go to Step 1.

\end{enumerate}
}
\end{algo}

Let $\nu$ and $Q$ be the mean and covariance matrix of the $n$ risky assets, respectively.
The  limited diversified mean-variance portfolio selection problem  can be formulated as
 \begin{align*}
{\rm(MV)}~~~~
\min~~
&x^TQx\\
{\rm s.t.}~~~
&\|x\|_0\le K,\\
&x\in X,
\end{align*}
where $\|x\|_0\le K$ is the cardinality constraint and
\[
X\triangleq\left\{x\in\Re^n~\middle|~ \nu^Tx\ge \rho,~\sum^n_{i=1}x_i=1,~0\le x_i\le u_i,~i=1,\ldots,n\right\},
\]
representing the constraints
of minimum return level, budget constraint and lower and upper bounds for $x_i$, respectively.

We consider the following equivalent reformulation of (MV):
 \begin{align*}
{\rm(RM)}~~~~
\min~~
&x^TQx\\
{\rm s.t.}~~~
&\sum_{i=1}^ny_i\leq K,\\
&0\leq y_{i}\leq 1,\\
&\nu^Tx\ge \rho,\\
&e^Tx-x^Ty\leq 0,\\
&\sum^n_{i=1}x_i=1,\\
&0\le x_i\le u_i,~i=1,\ldots,n,
\end{align*}
whose convex subproblem is
\begin{align*}
{\rm(MV_\mu(\bar{x},\bar{y}))}~~~~
\min~~
&x^TQx+\mu[(e-\bar{y},-\bar{x})(x,y)^T+\sum_{i=1}^{n}\bar{x}_i\bar{y}_i]\\
{\rm s.t.}~~~
&\sum_{i=1}^ny_i\leq K,\\
&0\leq y_{i}\leq 1,\\
&\nu^Tx\ge \rho,\\
&\sum^n_{i=1}x_i=1,\\
&0\le x_i\le u_i,~i=1,\ldots,n.
\end{align*}
\hspace*{6mm}Our test problems for (MV) consist of 90 instances,
where  parameters $Q$, $\nu$, $\rho$ and $u_i$ are created in the same way as in \cite{frangioni06,frangioni07}.
 There are  30 instances  for each
$n=200$, $300$ and $400$. The matrix $Q$ in the 30 instances for each problem size are generated
with different degrees of diagonal dominance. The parameters $\rho$  and $u_i$
 are uniformly drawn at random from intervals $[0.002, 0.01]$ and
$[0.375, 0.425]$, respectively. The  data files of these instances are available at:
{\tt http://www.di.unipi.it/optimize/Data/MV.html}.\\
\hspace*{6mm}Algorithm 1 is coded in Matlab (version R2012a) and executed on
  a PC equipped with Intel(R) Core(TM) i5-2520M CPU (2.50 GHz)
 and 4 GB of RAM.
  All the convex quadratic subproblems
in Algorithm 1 are solved by the QP solver in
{\tt CPLEX 12.3}  via the Matlab interface (see \cite{cplex}).

In our implementation, the initial  solution $x^0$ in Step 1 of Algorithm 1 is
obtained by solving the following convex quadratic programming:
\begin{align*}
  \min~~
  &x^TQx\\
  {\rm s.t.}~~~
  &\nu^Tx\ge \rho,\\
  &\sum^n_{i=1}x_i=1,\\
  &0\le x_i\le u_i,~i=1,\ldots,n.
\end{align*}
The initial vale of $\mu_0$ and $\varrho$ in Step 3 is set at $\mu_0=10$ and $\varrho=10$.  The value of stopping parameter $\epsilon$ in Step 3 is set at $\epsilon=10^{-7}$.
We set $p=\frac{1}{2}$ and $p=0.01$ for the $\ell_p$ function and the exponential function respectively, the same as in \cite{chen2010a,mangasarian1999}. \\
\hspace*{6mm}For $n=200$, $300$, and $400$, we solve problem $({\rm MV})$ by Algorithm 1.
From the numerical results as shown in Table \ref{tab1},  we could  get a local optimal
solution $(x^*,y^*)$ with $\|x^*\|_0= K$ and $y^*_i\in\{0,1\}$ ($i=1,\ldots,n$)  in a very short time, where $K$ is the cardinality number we set before the test.  This means that the penalty term in the  objective function of the  convex subproblem ${\rm AP}_\mu{({x}^k,{y}^k)}$ has become zero in the final solution of Algorithm 1.
 The average cardinality value is smaller than the cardinality value $K$ we set in some cases. It is because that in these cases there exist some problems that $\|x_0\|<K$, where $x_0$ is the initial solution we get in Step 1 in Algorithm 1. Comparing with ``SCA-$\ell_p$'' and ``SCA-exp'', which create different cardinality   solutions using different penalty parameters \cite{chen2010a,mangasarian1999}, our algorithm can set the cardinality number initially, i.e., we can select a desirable cardinality number $K$ before our test. \\
\hspace*{6mm}We use  ``SCA-$\ell_p$'' and ``SCA-exp'' methods to create 13 different sparse solutions respectively for $n=200$, 300, and 400.
 Using the cardinality of the sparse solution created by ``SCA-$\ell_p$'' and ``SCA-exp'' through different
 penalty parameters, we test the objective value, computing time and other items with the same cardinality of the sparse solution got by our algorithm. So the average cardinality in a line are equal to each other.
 The results show that when solving the same model, the time used by Algorithm 1 is much shorter than ``SCA-$\ell_p$'' and ``SCA-exp''.  The results are shown in Tables  \ref{tab2}, \ref{tab3}, \ref{tab4}, \ref{tab5}, \ref{tab6} and \ref{tab7}, where
the notations used are defined as follows:
\begin{itemize}
\item ``SCA-AP'', ``SCA-$\ell_p$'' and ``SCA-exp''  stand for Algorithm 1, the ``$\ell_p$'' successive  approximation  method   and the exponential approximation method respectively;
\item `` ${\rm K_{aver}}$ '' is the average value of cardinality (sparsity), which is the value of $\|x\|_0$,
of the sparse solutions generated by Algorithm 1 for   30 instances;
\item ``obj'' denotes the average objective function values $f(x)$ of
the sparse solutions generated by Algorithm 1 for
  30 instances;
\item  ``${\rm iter_a}$'' and ``${\rm iter_s}$'' denote the average number of iterations of Algorithm 1
and the average number of inner iterations in solving the subproblem
for the 30 instances, respectively;
\item ``${\rm time_a}$'' and ``${\rm time_s}$''  denote the average CPU time of Algorithm 1
and the average computing time in solving the subproblem at each iteration
for the 30 instances, respectively.
\end{itemize}

Compared with the regularized method, our method can get a local optimal solution which exactly satisfies the cardinality constraint.
The average CPU time of Algorithm 1 are much  shorter than those of ``SCA-$\ell_p$'' and ``SCA-exp''.
In our experiments, we find that the average CPU time of Algorithm 1 is much  shorter than  ``SCA-$\ell_p$'' and ``SCA-exp'',
and the average objective values are much better  than  ``SCA-$\ell_p$'' and ``SCA-exp'' when $K$ is big.
In the case where $K$ is small, the average objective values  from Algorithm 1 are little bit larger than these of ``SCA-$\ell_p$'' and ``SCA-exp''.
The difference of average objective values between Algorithm 1 and ``SCA-$\ell_p$'' or ``SCA-exp'' is no more than 2. All the above results show the effectiveness of
our algorithm.

\section{Conclusions}
We have presented some prominent properties of the cardinality constrained optimization program under some conditions. In particular, we have developed an equivalent reformulation for the optimization  problem with a sparsity constraint and nonnegative variables. Based on this reformulation,  we have further constructed a successive convex approximation (SCA) method and established the convergence of
the sequence of approximate solutions to a KKT point of
the original problem. We finally confirmed the effectiveness of our algorithm from the computational results of the limited diversified mean-variance
portfolio selection problem in our numerical tests.


\begin{table}[t!] \vspace*{0.0in} \centering
\caption{\label{tab1}
Numerical results for SCA-RP method}
{\small \tabcolsep=1.2pt
\begin{tabular}{ccccccccccccc}
\hline \noalign{\smallskip}
\multicolumn{2}{c}{ } &&
    \multicolumn{6}{c}{SCA-RP }  \\
      \noalign{\smallskip}
  \cline{3-13}
 \noalign{\smallskip}
 $n$ &&$K$ &&${\rm K_{aver}}$&&$\mu_0$&& obj &${\rm iter_a}$ &${\rm time_a}$ &${\rm iter_s}$ & ${\rm time_s}$   \\

 \hline  \noalign{\smallskip}
 200 && 150 &&  139.93 && 10&&      17.85 &     4 & 0.62 &    18 & 0.18 \\
 200 && 130 &&  123.83 && 10&&      18.66 &     2 & 0.31 &    12 & 0.15 \\
 200 && 100 &&  98.30 &&  10&&      20.55 &    4 & 0.63 &    12 & 0.15  \\
 200 &&80 &&    79.40 &&  10&&      22.84 &    2 & 0.32 &    13 & 0.16  \\
 200 &&60 &&    60.00 &&  10&&      27.25 &    2 & 0.32 &    12 & 0.16  \\
 200 &&40 &&    40.00 &&  10&&      36.70 &    3 & 0.43 &    12 & 0.16  \\
 200 &&20 &&    20.00 &&  10&&      64.53 &    3 & 0.41 &    14 & 0.16  \\
 200 &&15 &&    15.00 &&  10&&      82.69 &    11 & 1.88 &    15 & 0.16 \\
 200&& 10 &&    10.00 &&  10&&      119.07 &    8 & 1.56 &    17 & 0.18 \\
 200&&5 &&       5.00 &&  10&&      227.71 &    8 & 1.41 &    18 & 0.18 \\
 \hline \noalign{\smallskip}
    300 && 200&&   180.93 && 10&& 19.55 &        2 & 0.65 &    16 & 0.26 \\
    300 && 150&&   139.57 && 10&& 21.93 &        2 & 0.39 &    11 & 0.21 \\
    300 &&100 &&   95.37 &&  10&& 27.14 &       2 & 0.42 &    12 & 0.22  \\
    300 &&80 &&    77.37 &&  10&& 31.12 &       2 & 0.45 &    12 & 0.23  \\
    300 &&60 &&    58.93 &&  10&& 37.83 &       3 & 0.60 &    12 & 0.23  \\
    300 && 40 &&   39.70 &&  10&& 51.41 &       3 & 0.64 &    13 & 0.24  \\
    300 &&20&&     20.00 &&  10&& 92.90 &      15 & 4.06 &    16 & 0.26  \\
    300 &&15 &&    15.00 &&  10&& 120.69 &     8 & 2.21 &    16 & 0.26   \\
    300 && 10&&    10.00 &&  10&& 175.58 &     6 & 1.59 &    18 & 0.27   \\
    300 && 5&&      5.00 &&  10&& 339.63 &     7 & 1.82 &    16 & 0.27   \\
 \hline \noalign{\smallskip}
 400&&300&&  264.50 && 10&& 23.32 &      2 & 0.79 &    18 & 0.43 \\
 400&&250&&  227.83 && 10&& 24.35 &      2 & 0.69 &    13 & 0.36 \\
 400&&200&&  189.73 && 10&& 26.00 &      2 & 0.72 &    13 & 0.38 \\
 400&&150&&  145.10 && 10&& 28.84 &      2 & 0.78 &    13 & 0.39 \\
 400&&100&&  98.37 &&  10&& 35.25 &     4 & 1.78 &    14 & 0.43  \\
 400&&80&&   79.03 &&  10&& 40.41 &     4 & 1.80 &    13 & 0.40  \\
 400&&60&&   59.70 &&  10&& 49.32 &     3 & 1.27 &    14 & 0.47  \\
 400 &&40&&  40.00 &&  10&& 67.00 &     3 & 1.40 &    16 & 0.50  \\
 400 &&20&&  20.00 &&  10&& 123.66 &    14 & 8.13 &    17 & 0.53  \\
 400 && 15&& 15.00 &&  10&& 160.53 &     8 & 4.29 &    17 & 0.52  \\
 400 &&10 && 10.00 &&   10&&233.74 &    12 & 6.99 &    19 & 0.55  \\
 400 && 5&&   5.00 &&   10&& 453.40 &    7 & 3.84 &    18 & 0.54  \\
 \hline \noalign{\smallskip}
\end{tabular}
}
\end{table}

\begin{table}[t!] \vspace*{0.0in} \centering
\caption{\label{tab2}
Comparison between SCA-RP and SCA-exp with $n=200$}
{\small \tabcolsep=1.2pt
\begin{tabular}{cccccccccccccccccc}
\hline \noalign{\smallskip}
\multicolumn{6}{c}{SCA-RP } &&
 \multicolumn{6}{c}{SCA-exp}  \\
      \noalign{\smallskip}
 \cline{1-6} \cline{8-13}
 \noalign{\smallskip}
 ${\rm K_{aver}}$ &obj & ${\rm iter_a}$ &${\rm time_a}$ &${\rm iter_s}$ & ${\rm time_s}$ &&
  ${\rm K_{aver}}$  &obj & ${\rm iter_a}$ &${\rm time_a}$ &${\rm iter_s}$ & ${\rm time_s}$ \\

 \hline  \noalign{\smallskip}

 30.23 & 42.43 &     3 & 0.38 &    14 & 0.13   && 30.23 & 41.64 &    30 & 2.81 &    10 & 0.09  \\
 40.33 & 33.96 &     2 & 0.31 &    13 & 0.16    && 40.33 & 33.59 &    49 & 4.93 &    10 & 0.10  \\
 50.17 & 29.25 &     2 & 0.23 &    13 & 0.12   && 50.17 & 29.23 &    58 & 5.96 &    10 & 0.11  \\
 60.20 & 26.57 &     2 & 0.22 &    12 & 0.11  && 60.20 & 26.82 &    84 &10.11 &    10 & 0.12  \\
 70.43 & 25.18 &     2 & 0.26 &    13 & 0.13  && 70.43 & 25.77 &    86 &10.50 &    11 & 0.12  \\
 80.30 & 23.89 &     2 & 0.26 &    11 & 0.14  && 80.30 & 24.14 &    49 & 5.89 &    10 & 0.12  \\
 88.93 & 21.80 &     2 & 0.21 &    11 & 0.12  && 88.93 & 21.74 &    55 & 6.74 &    10 & 0.12  \\
 98.10 & 20.31 &     2 & 0.20 &    12 & 0.11  && 98.10 & 20.27 &    68 & 8.02 &    10 & 0.12  \\
 110.03 & 19.06 &     2 & 0.20 &    12 & 0.11  && 110.03 & 19.17 &    84 & 7.43 &    11 & 0.09  \\
 121.13 & 18.36 &     2 & 0.21 &    12 & 0.11 & & 121.13 & 18.63 &   106 & 9.49 &    11 & 0.09  \\
 130.80 & 17.94 &     2 & 0.24 &    14 & 0.13 & & 130.80 & 18.29 &   104 & 9.26 &    11 & 0.09  \\
 140.47 & 17.59 &     2 & 0.22 &    13 & 0.12 & & 140.47 & 17.96 &    93 & 8.14 &    11 & 0.09  \\
 149.60 & 17.25 &     2 & 0.17 &    11 & 0.10  & & 149.60 & 17.53 &    72 & 6.28 &    10 & 0.08  \\

\hline \noalign{\smallskip}
\end{tabular}
}
\end{table}

\begin{table}[t!] \vspace*{0.0in} \centering
\caption{\label{tab3}
Comparison between SCA-RP and SCA-$\ell_p$  with $n=200$}
{\small \tabcolsep=1.2pt
\begin{tabular}{cccccccccccccccccc}
\hline \noalign{\smallskip}
\multicolumn{6}{c}{SCA-RP } &&
 \multicolumn{6}{c}{SCA-$\ell_p$}  \\
      \noalign{\smallskip}
 \cline{1-6} \cline{8-13}
 \noalign{\smallskip}
  ${\rm K_{aver}}$  &obj & ${\rm iter_a}$ &${\rm time_a}$ &${\rm iter_s}$ & ${\rm time_s}$ &&
  ${\rm K_{aver}}$  &obj & ${\rm iter_a}$ &${\rm time_a}$ &${\rm iter_s}$ & ${\rm time_s}$ \\

 \hline  \noalign{\smallskip}

 30.43 & 42.47 &     2 & 0.37 &    14 & 0.16  & & 30.43 & 41.75 &    54 & 4.04 &    10 & 0.07   \\
 40.50 & 33.95 &     2 & 0.32 &    13 & 0.16  &         & 40.50 & 33.53 &    59 & 4.83 &    10 & 0.08   \\
 49.67 & 29.23 &     2 & 0.31 &    13 & 0.16  &          & 49.67 & 29.08 &    49 & 4.04 &    10 & 0.08   \\
 59.43 & 25.92 &     2 & 0.32 &    12 & 0.16  &          & 59.43 & 25.82 &    49 & 3.90 &     9 & 0.08   \\
 70.50 & 23.37 &     2 & 0.30 &    12 & 0.15  &         & 70.50 & 23.32 &    50 & 3.93 &     9 & 0.08   \\
 80.93 & 21.73 &     2 & 0.32 &    13 & 0.16  &          & 80.93 & 21.69 &    59 & 4.62 &     9 & 0.08   \\
 91.20 & 20.50 &     2 & 0.32 &    14 & 0.16  &   & 91.20 & 20.47 &    56 & 4.37 &     9 & 0.08   \\
 100.90 & 19.59 &     2 & 0.31 &    13 & 0.16  &         & 100.90 & 19.58 &    55 & 4.29 &     9 & 0.08   \\
 109.90 & 18.89 &     2 & 0.31 &    12 & 0.15  &         & 109.90 & 18.89 &    58 & 4.55 &     9 & 0.08   \\
 119.77 & 18.22 &     2 & 0.27 &    12 & 0.15  &         & 119.77 & 18.23 &    49 & 3.92 &     9 & 0.08   \\
 129.57 & 17.73 &     2 & 0.27 &    11 & 0.15  &         & 129.57 & 17.75 &    45 & 3.57 &     9 & 0.08   \\
 140.83 & 17.31 &     2 & 0.30 &    14 & 0.16  &         & 140.83 & 17.35 &    42 & 3.32 &     9 & 0.08   \\
 149.10 & 17.07 &     2 & 0.28 &    13 & 0.15    && 149.10 & 17.11 &    47 & 3.76 &     9 & 0.08   \\

\hline \noalign{\smallskip}
\end{tabular}
}
\end{table}

\begin{table}[t!] \vspace*{0.0in} \centering
\caption{\label{tab4}
Comparison between SCA-RP and SCA-exp with $n=300$}
{\small \tabcolsep=1.2pt
\begin{tabular}{cccccccccccccccccc}
\hline \noalign{\smallskip}
\multicolumn{6}{c}{SCA-RP } &&
 \multicolumn{6}{c}{SCA-exp}  \\
      \noalign{\smallskip}
 \cline{1-6} \cline{8-13}
 \noalign{\smallskip}
  ${\rm K_{aver}}$  &obj & ${\rm iter_a}$ &${\rm time_a}$ &${\rm iter_s}$ & ${\rm time_s}$ &&
  ${\rm K_{aver}}$  &obj & ${\rm iter_a}$ &${\rm time_a}$ &${\rm iter_s}$ & ${\rm time_s}$ \\

 \hline  \noalign{\smallskip}
 30.30 & 60.15 &     4 & 1.01 &    15 & 0.22 &  & 30.30 & 59.31 &    35 & 5.55 &    10 & 0.16 \\
 40.03 & 47.94 &     3 & 0.56 &    14 & 0.20 &   & 40.03 & 47.08 &    48 & 8.18 &    10 & 0.18 \\
 50.50 & 40.38 &     2 & 0.41 &    14 & 0.20 &   & 50.50 & 39.73 &    64 &12.42 &    10 & 0.19 \\
 60.53 & 35.91 &     2 & 0.39 &    13 & 0.20 &   & 60.53 & 35.62 &    86 &14.02 &    10 & 0.17 \\
 70.30 & 33.07 &     2 & 0.39 &    13 & 0.19 &  & 70.30 & 33.20 &    93 &17.69 &    10 & 0.20 \\
 79.77 & 31.32 &     2 & 0.39 &    13 & 0.20 &   & 79.77 & 31.89 &    92 &19.05 &    11 & 0.20 \\
 88.40 & 30.25 &     2 & 0.41 &    14 & 0.21 &   & 88.40 & 31.23 &   103 &20.45 &    11 & 0.20 \\
 97.77 & 29.48 &     2 & 0.41 &    14 & 0.20 &   & 97.77 & 30.78 &    90 &18.55 &    11 & 0.20 \\
 107.93 & 28.67 &     2 & 0.43 &    15 & 0.22 &  & 107.93 & 30.18 &    92 &17.06 &    11 & 0.18 \\
 120.33 & 27.69 &     2 & 0.42 &    13 & 0.23 &  & 120.33 & 28.90 &    63 &11.00 &    10 & 0.17 \\
 130.77 & 25.23 &     2 & 0.38 &    12 & 0.22 &  & 130.77 & 25.39 &    54 & 9.61 &    10 & 0.18 \\
 140.43 & 22.93 &     2 & 0.37 &    12 & 0.22 &  & 140.43 & 22.99 &    69 &13.65 &    10 & 0.19 \\
 149.40 & 21.49 &     2 & 0.38 &    12 & 0.22 &  & 149.40 & 21.62 &    87 &17.23 &    10 & 0.20 \\

\hline \noalign{\smallskip}
\end{tabular}
}
\end{table}

\begin{table}[t!] \vspace*{0.0in} \centering
\caption{\label{tab5}
Comparison between SCA-RP and SCA-$\ell_p$ with $n=300$}
{\small \tabcolsep=1.2pt
\begin{tabular}{cccccccccccccccccc}
\hline \noalign{\smallskip}
\multicolumn{6}{c}{SCA-RP } &&
 \multicolumn{6}{c}{SCA-$\ell_p$}  \\
      \noalign{\smallskip}
 \cline{1-6} \cline{8-13}
 \noalign{\smallskip}
  ${\rm K_{aver}}$  &obj & ${\rm iter_a}$ &${\rm time_a}$ &${\rm iter_s}$ & ${\rm time_s}$ &&
  ${\rm K_{aver}}$  &obj & ${\rm iter_a}$ &${\rm time_a}$ &${\rm iter_s}$ & ${\rm time_s}$ \\

 \hline  \noalign{\smallskip}

30.43 & 60.73 &     8 & 2.21 &    16 & 0.25 &   & 30.43 & 59.85 &    56 & 8.21 &    11 & 0.15 \\
40.87 & 47.44 &     3 & 0.66 &    14 & 0.24  &   & 40.87 & 46.69 &    60 & 9.83 &    10 & 0.16 \\
49.27 & 40.98 &     2 & 0.49 &    14 & 0.23  &   & 49.27 & 40.32 &    56 & 8.96 &    11 & 0.16 \\
60.63 & 35.00 &     2 & 0.45 &    13 & 0.22 &    & 60.63 & 34.51 &    56 & 8.98 &    10 & 0.16 \\
69.83 & 31.62 &     2 & 0.44 &    13 & 0.22 &    & 69.83 & 31.29 &    57 & 8.92 &    10 & 0.16 \\
80.37 & 28.71 &     2 & 0.44 &    13 & 0.22 &   & 80.37 & 28.52 &    60 & 9.84 &    10 & 0.17 \\
89.17 & 26.83 &     2 & 0.43 &    12 & 0.21  &   & 89.17 & 26.70 &    59 & 9.48 &    10 & 0.16 \\
99.53 & 25.15 &     2 & 0.43 &    12 & 0.21  &   & 99.53 & 25.05 &    58 & 9.36 &    10 & 0.16 \\
109.70 & 23.78 &     2 & 0.45 &    12 & 0.22  &  & 109.70 & 23.74 &    55 & 8.79 &     9 & 0.16 \\
119.57 & 22.71 &     2 & 0.46 &    13 & 0.23  &  & 119.57 & 22.70 &    56 & 8.78 &     9 & 0.16 \\
130.90 & 21.67 &     2 & 0.48 &    15 & 0.24  &  & 130.90 & 21.70 &    64 &10.05 &     9 & 0.16 \\
140.03 & 20.97 &     2 & 0.46 &    13 & 0.22   &         & 140.03 & 21.02 &    61 & 9.74 &     9 & 0.16 \\
149.77 & 20.33 &     2 & 0.46 &    13 & 0.23  &  & 149.77 & 20.40 &    62 & 9.78 &     9 & 0.16 \\
\hline \noalign{\smallskip}
\end{tabular}
}
\end{table}

\begin{table}[t!] \vspace*{0.0in} \centering
\caption{\label{tab6}
Comparison between SCA-RP and SCA-exp with $n=400$}
{\small \tabcolsep=1.2pt
\begin{tabular}{cccccccccccccccccc}
\hline \noalign{\smallskip}
\multicolumn{6}{c}{SCA-RP } &&
 \multicolumn{6}{c}{SCA-exp}  \\
      \noalign{\smallskip}
 \cline{1-6} \cline{8-13}
 \noalign{\smallskip}
  ${\rm K_{aver}}$  &obj & ${\rm iter_a}$ &${\rm time_a}$ &${\rm iter_s}$ & ${\rm time_s}$ &&
  ${\rm K_{aver}}$  &obj & ${\rm iter_a}$ &${\rm time_a}$ &${\rm iter_s}$ & ${\rm time_s}$ \\

 \hline  \noalign{\smallskip}

 30.13 & 79.02 &    11 & 5.01 &    17 & 0.45 & &     30.13 & 77.82 &    29 & 8.03 &    11 & 0.28 \\
  40.30 & 62.22 &    6 & 3.15 &    17 & 0.47 & &    40.30 & 61.28 &    45 &13.36 &    10 & 0.30 \\
  50.83 & 52.61 &    3 & 1.49 &    15 & 0.48 & &    50.83 & 51.99 &    64 &18.60 &    10 &0.29 \\
  62.37 & 46.60 &    2 & 0.96 &    15 & 0.48 & &    62.37 & 46.21 &    80 &23.79 &    10 &  0.29 \\
  70.30 & 44.15 &    2 & 0.97 &    15 & 0.48 & & 70.30 & 44.06 &    87 &25.88 &   10 & 0.30 \\
  80.30 & 42.04 &    2 & 0.91 &    14 & 0.46 & &    80.30 & 42.30 &    89 &28.02 &    11  & 0.31 \\
  95.63 & 40.35 &    2 & 0.92 &    15 & 0.46 & &    95.63 & 41.32 &    87 &27.24 &    11 & 0.30 \\
  113.30 & 39.19 &    2 & 0.98 &    16 & 0.49& &    113.30 & 40.54 &    84 &25.93 &    11 & 0.30 \\
  128.53 & 38.07 &    2 & 0.94 &    15 & 0.48& &   128.53 & 39.17 &    72 &21.70 &    11 & 0.30 \\
  132.97 & 37.05 &    2 & 0.83 &    14 & 0.45& &   132.97 & 37.33 &    59 &17.43 &    10 & 0.29 \\
  143.50 & 33.45 &    2 & 0.81 &    13 & 0.44& &   143.50 & 33.04 &    59 &17.64 &    10 & 0.30 \\
  150.50 & 31.39 &    2 & 0.77 &    13 & 0.42& &   150.50 & 31.01 &    72 &21.97 &    11 & 0.30 \\
  163.20 & 28.90 &    2 & 0.77 &    13 & 0.42& &   163.20 & 28.72 &    79 &26.04 &    11 & 0.32\\

 \hline \noalign{\smallskip}
\end{tabular}
}
\end{table}

\begin{table}[t!] \vspace*{0.0in} \centering
\caption{\label{tab7}
Comparison between SCA-RP and SCA-$\ell_p$ with $n=400$}
{\small \tabcolsep=1.2pt
\begin{tabular}{cccccccccccccccccc}
\hline \noalign{\smallskip}
\multicolumn{6}{c}{SCA-RP } &&
 \multicolumn{6}{c}{SCA-$\ell_p$}  \\
      \noalign{\smallskip}
 \cline{1-6} \cline{8-13}
 \noalign{\smallskip}
  ${\rm K_{aver}}$  &obj & ${\rm iter_a}$ &${\rm time_a}$ &${\rm iter_s}$ & ${\rm time_s}$ &&
  ${\rm K_{aver}}$  &obj & ${\rm iter_a}$ &${\rm time_a}$ &${\rm iter_s}$ & ${\rm time_s}$ \\

 \hline  \noalign{\smallskip}
  60.87 & 45.89 &    2 & 0.88 &    15 & 0.37  & & 60.87 & 45.33 &    69 &22.27 &    10 & 0.32 \\
  70.10 & 41.71 &    2 & 0.76 &    15 & 0.38  & & 70.10 & 41.19 &    67 &21.96 &    10 & 0.33 \\
  80.17 & 38.27 &    2 & 0.71 &    14 & 0.36  &  & 80.17 & 37.81 &    59 &19.55 &    10 & 0.33 \\
  90.30 & 35.63 &    2 & 0.74 &    14 & 0.37  & & 90.30 & 35.17 &    59 &19.16 &    10 & 0.33 \\
  101.87 & 33.30 &    2 & 0.71 &    13 & 0.35  & & 101.87 & 32.89 &    57 &19.05 &    10 & 0.33 \\
  110.77 & 31.84 &    2 & 0.65 &    13 & 0.32  & & 110.77 & 31.49 &    64 &21.95 &    10 & 0.34 \\
  121.40 & 30.36 &    2 & 0.64 &    12 & 0.32  &  & 121.40 & 30.08 &    62 &20.09 &    10 & 0.33 \\
  126.57 & 29.75 &    2 & 0.64 &    13 & 0.32  &  & 126.57 & 29.49 &    62 &20.42 &    10 & 0.33 \\
  136.10 & 28.75 &    2 & 0.66 &    12 & 0.33  &  & 136.10 & 28.55 &    62 &19.82 &    10 & 0.32 \\
  140.93 & 28.29 &    2 & 0.69 &    13 & 0.35   &  & 140.93 & 28.12 &    57 &18.25 &    10 & 0.32 \\
  150.50 & 27.49 &    2 & 0.74 &    15 & 0.37  & & 150.50 & 27.35 &    66 &21.44 &     9 & 0.32 \\
  164.47 & 26.49 &    2 & 0.73 &    14 & 0.36  &  & 164.47 & 26.39 &    57 &18.27 &    10 & 0.32 \\
  182.43 & 25.46 &    2 & 0.72 &    14 & 0.36   &  & 182.43 & 25.40 &    67 &21.51 &     9 & 0.32 \\

  \hline \noalign{\smallskip}
\end{tabular}
}
\end{table}

\end{document}